\UseRawInputEncoding 
\documentclass[final]{siamltex}%
\usepackage{eurosym}
\usepackage{algorithm}
\usepackage{algpseudocode}
\usepackage{amsfonts}
\usepackage{amsmath}
\usepackage{amssymb}
\usepackage{comment}
\usepackage{enumerate}
\usepackage{graphicx}
\usepackage{mathtools}
\usepackage{multirow}
\usepackage{url}
\usepackage{verbatim}
\usepackage[pdftex]{thumbpdf}
\usepackage[pdftex,
pdfstartview=FitH,bookmarks=true,bookmarksnumbered=true,bookmarksopen=true,hypertexnames=false,breaklinks=true,
colorlinks=true,  linkcolor=blue,anchorcolor=blue,
citecolor=blue,filecolor=blue,  menucolor=blue]{hyperref}
\usepackage{comment}%
\providecommand{\U}[1]{\protect\rule{.1in}{.1in}}

\providecommand{\U}[1]{\protect\rule{.1in}{.1in}}

\newtheorem{remark}[theorem]{Remark}

\begin{document}

\title{A stochastic Galerkin method \\with adaptive time-stepping for the \\Navier\textendash Stokes equations\thanks{This work was supported by the
U.\thinspace\ S.\thinspace\ National Science Foundation under grant
DMS1913201. Most of the computations were performed at the Center for
Computational Mathematics, University of Colorado Denver, and the
support of Prof. Jan Mandel is greatly appreciated.
Part of the work was completed while Randy Price was student at the University of Maryland, Baltimore County.}}
\author{Bed\v{r}ich Soused\'{\i}k\thanks{Department of Mathematics and Statistics,
University of Maryland, Baltimore County, 1000 Hilltop Circle, Baltimore,
MD~21250 (\texttt{sousedik@umbc.edu}).}
\and Randy Price\thanks{The Center for Mathematics and Artificial Intelligence and the Center for Computational Fluid Dynamics,
George Mason University, Fairfax, VA~22030
(\texttt{rprice25@gmu.edu})}.}
\maketitle

\begin{abstract}
We study the time-dependent Navier\textendash Stokes equations in the context
of stochastic finite element discretizations. Specifically, we assume that the
viscosity is a random field given in the form of a generalized polynomial
chaos expansion, and we use the stochastic Galerkin method to extend the
methodology from [D. A. Kay et al., \textit{SIAM J. Sci. Comput.} 32(1), pp.
111--128, 2010] into this framework. For the resulting stochastic problem, we
explore the properties of the resulting stochastic solutions, and we also
compare the results with that of Monte Carlo and stochastic collocation. Since
the time-stepping scheme is fully implicit, we also propose strategies for
efficient solution of the stochastic Galerkin linear systems using a
preconditioned Krylov subspace method. The effectiveness of the stochastic
Galerkin method is illustrated by numerical experiments.

\end{abstract}



\pagestyle{myheadings} \thispagestyle{plain} \markboth{B.\ SOUSED\'{I}K AND R. PRICE}{STOCHASTIC GALERKIN METHOD FOR NAVIER\textendash STOKES EQUATIONS}

\section{Introduction}

\label{sec:introduction}Models of mathematical physics are commonly based on
partial differential equations (PDEs). In this study, we focus on the most
popular PDE model in fluid mechanics, which is the Navier--Stokes
equation~\cite{Elman-2014-FEF,LeMaitre-2010-SMU}. We consider a stochastic
version of the model: we assume that the viscosity is given by a generalized
polynomial chaos (gPC) expansion, we discretize the problem using spectral
stochastic finite elements see,
e.g.,~\cite{Ghanem-1991-SFE,LeMaitre-2010-SMU,Lord-2014-ICS,Xiu-2010-NMS}, and
we wish to find the gPC\ expansion of the solution. The steady-state version
of this problem was studied
in~\cite{Lee-2019-LRS,Powell-2012-PSS,Sousedik-2016-SGM}, and our focus here
is on the time-dependent counterpart. Our approach to time discretization is
built on the fully implicit scheme with adaptive time-stepping strategy, which
was developed for the deterministic Navier--Stokes equation by Kay et
al.~\cite{Kay-2010-ATS}, see also~\cite{Gresho-2008-ATS}. We extend their
scheme in the stochastic Galerkin framework, and in particular we show that
the physics inspired time-stepping strategy can be also adapted to this
framework. The scheme is fully implicit, and so each time step entails a solve
with the stochastic Galerkin matrix. This typically leads to very large
systems of linear equations, for which use of direct solvers may be
prohibitive, and therefore the method could potentially be quite
computationally expensive. There are other approaches to time stepping see,
e.g.,~\cite{Almgren-1996-NMI,Elman-2014-FEF,Karniadakis-1991-HSM} which may
appear more appealing. Nevertheless, finally we also show that the iterative
solvers by Soused\'{\i}k and Elman~\cite{Sousedik-2016-SGM}, which are based
on preconditioned Krylov subspace methods, are quite effective for the
implicit time discretizations of the time-dependent Navier--Stokes problem as
well.

Some aspects of the gPC\ methods for time-dependent problems were studied in
literature see, e.g.,~\cite{Gunzburger-2019-EFR,Xiu-2003-MUF,Xiu-2003-NSA}. In
particular, long-term integration was addressed by Gerritsma et
al.~\cite{Gerritsma-2010-TDG}, Heuveline, Schick and
Song~\cite{Heuveline-2014-HGP,Schick-2011-UQS,Song-2017-MPP},
Wilkins~\cite{Wilkins-2016-ECE}, \"{O}zen nad
Bal~\cite{Ozen-2016-DPC,Ozen-2017-DPC},\ and most recently by Esquivel et
al.~\cite{Esquivel-2021-FDS}, among others. Methods for flows exhibiting
uncertain periodic dynamics were proposed, e.g., by Bonnaire et
al.~\cite{Bonnaire-2021-OGP}, Lacour et al.~\cite{Lacour-2021-DSF} and Schick
et al.~\cite{Schick-2016-NGM}. These methods typically entail time-dependent
or other variants of gPC\ expansions that are tailored to the changing
character of the solution.
Nevertheless, here we use a time-independent gPC basis, which turns out to be
sufficient for the transient problems considered in our numerical experiments.
Therefore, all these techniques can be viewed as complementary to the present
study. We also note that Elman and Su~\cite{Elman-2020-LRS} proposed a
low-rank stochastic Galerkin solver based on monolithic (all-at-once) time
discretization of the Navier--Stokes problem, however their scheme is based on
a constant timestep.

Finally, we remark on possible interpretations of the Navier--Stokes problem
with stochastic viscosity. In such case, the Reynolds number defined as
\[
\operatorname*{Re}(\xi)=\frac{UL}{\nu(\xi)},
\]
where $\nu>0$ is the viscosity, $U$ is the characteristic velocity and $L$ is
the characteristic length, is also stochastic. The possible interpretations of
such setup are discussed by Powell and Silvester in~\cite{Powell-2012-PSS}:
for example, assuming fixed geometry, the stochastic viscosity is equivalent
to Reynolds number being stochastic, which may correspond to a scenario when
the volume of fluid moving into the channel is uncertain.

The paper is organized as follows. In Section~\ref{sec:model} we recall the
algorithm for the deterministic problem, in Section~\ref{sec:stochastic}\ we
formulate the algorithm for the stochastic problem using both the stochastic
Galerkin and sampling methods, in Section~\ref{sec:numerical} we report
results of numerical experiments and provide details about the preconditioning
of the Oseen problem, and finally in Section~\ref{sec:conclusion} we summarize
and conclude our work.

\section{Algorithm for the deterministic problem}

\label{sec:model}We first recall the algorithm for the deterministic problem
following Kay et al.~\cite{Kay-2010-ATS}. Let $D\subset%
\mathbb{R}
^{2}$ be a physical domain, and let $T>0$ denote a stopping time. We wish to
solve the time-dependent Navier--Stokes equation in $D\times\lbrack0,T]$,
where$\ \left(  \vec{u},p\right)  $ denote the fluid velocity and pressure,
and $\nu\equiv\nu(x)>0$\footnote{The assumption $\nu(x)\neq
\operatorname*{const}$ is the only difference from the setup
in~\cite{Kay-2010-ATS} in this section.} is the viscosity parameter, written
as
\begin{align}
\frac{\partial\vec{u}}{\partial t} &  =f(\nu,\vec{u},p),\qquad f(\nu,\vec
{u},p)=\nu\nabla^{2}\vec{u}-\vec{u}\cdot\nabla\vec{u}-\nabla
p,\label{eq:NS1-a}\\
-\nabla\cdot\vec{u} &  =0,\label{eq:NS1-b}%
\end{align}
with boundary and initial conditions given on $\partial D=\overline{\Gamma
}_{D}\cup\overline{\Gamma}_{N}$ as
\begin{align}
\vec{u} &  =\vec{g},\quad\text{on }\Gamma_{D}\times\lbrack
0,T],\label{eq:bc-Dirichlet}\\
\nu\nabla\vec{u}\cdot\vec{n}-p\vec{n} &  =\vec{0},\quad\text{on }\Gamma
_{N}\times\lbrack0,T],\label{eq:bc-do-nothing}\\
\vec{u}(\vec{x},0) &  =\vec{u}_{0}(\vec{x}),\quad\text{in }D.\label{eq:ic}%
\end{align}
The initial velocity field is assumed to satisfy the incompressibility
constraint, that is $\nabla\cdot\vec{u}_{0}=0$. We also assume that
$\Gamma_{N}$ has nonzero measure so that the pressure is uniquely specified,
and to this end we will use the outflow (\emph{do-nothing}) boundary
condition. We begin by recalling the implicit trapezoid rule (TR) as
\[
\vec{u}_{t}\approx\frac{\vec{u}^{n+1}-\vec{u}^{n}}{k_{n+1}}=\frac{1}{2}\left[
f_{n+1}+f_{n}\right]  ,
\]
where $k_{n+1}=t_{n+1}-t_{n}$. Then $2\vec{u}_{t}\approx2\left(  \vec{u}%
^{n+1}-\vec{u}^{n}\right)  /k_{n+1}=f_{n+1}+f_{n}$ and (\ref{eq:NS1-a}%
)--(\ref{eq:NS1-b}) can be written as
\begin{align}
\frac{2}{k_{n+1}}\vec{u}^{n+1}-\nu\nabla^{2}\vec{u}^{n+1}+\vec{u}^{n+1}%
\cdot\nabla\vec{u}^{n+1}+\nabla p^{n+1} &  =\frac{2}{k_{n+1}}\vec{u}^{n}%
+\frac{\partial\vec{u}^{n}}{\partial t},\label{eq:NS2-a}\\
-\nabla\cdot\vec{u}^{n+1} &  =0.\label{eq:NS2-b}%
\end{align}
The nonlinear term is linearized as $\vec{u}^{n+1}\cdot\nabla\vec{u}%
^{n+1}\approx\vec{w}^{n+1}\cdot\nabla\vec{u}^{n+1}$. The linearization is
based on extrapolation $\left(  \vec{w}^{n+1}-\vec{u}^{n}\right)
/k_{n+1}=\left(  \vec{u}^{n}-\vec{u}^{n-1}\right)  /k_{n}$, from which we
find
\begin{equation}
\vec{w}^{n+1}=\left(  1+k_{n+1}/k_{n}\right)  \vec{u}^{n}-\left(
k_{n+1}/k_{n}\right)  \vec{u}^{n-1}.\label{eq:w}%
\end{equation}
Next, let~$(V_{D},Q_{D})$ denote a pair of spaces satisfying the inf-sup
condition and let $V_{E}$ be an extension of~$V_{D}$ containing velocity
vectors that satisfy the Dirichlet boundary
conditions~\cite{Brezzi-1991-MHF,Elman-2014-FEF,Girault-1986-FEM}. The mixed
variational formulation of (\ref{eq:NS2-a})--(\ref{eq:NS2-b}) is: find
$(\vec{u}^{n+1},p^{n+1})\in V_{E}\times Q_{D}$, for a given pair $(\vec{u}%
^{n},p^{n})$, such that
\begin{align}
\frac{2}{k_{n+1}}\int_{D}\vec{u}^{n+1}\vec{v}+\int_{D}\nu\,\nabla\vec{u}^{n+1}
&  :\nabla\vec{v}+\int_{D}\left(  \vec{w}^{n+1}\cdot\nabla\vec{u}%
^{n+1}\right)  \vec{v}-\int_{D}p^{n+1}\left(  \nabla\cdot\vec{v}\right)
\label{eq:NS3-a}\\
&  =\frac{2}{k_{n+1}}\int_{D}\vec{u}^{n}\vec{v}+\int_{D}\frac{\partial\vec
{u}^{n}}{\partial t}\vec{v},\nonumber\\
-\int_{D}q\left(  \nabla\cdot\vec{u}^{n+1}\right)   &  =0,\label{eq:NS3-b}%
\end{align}
for all $(\vec{v},q)\in V_{D}\times Q_{D}$. We note that $p^{n+1}$ is not
needed for subsequent time steps. Next, we recall the three ingredients of the
algorithm as discussed in \cite{Kay-2010-ATS}: time integration, time-step
selection and stabilization of the integrator.

\paragraph{Time integration}

Substituting $\vec{u}^{n+1}=\vec{u}^{n}+k_{n+1}\vec{d}^{n}$ into
(\ref{eq:NS3-a})--(\ref{eq:NS3-b}), rearranging and using $\int_{D}q\left(
\nabla\cdot\vec{u}^{n}\right)  =0$, we get the so-called \emph{discrete Oseen
problem}: given $\vec{u}^{n}$, $\partial\vec{u}^{n}/\partial t$\ and the
boundary update $\vec{g}\vcentcolon=(\vec{g}^{n+1}-\vec{g}^{n})/k_{n+1}$, we
first compute $\left(  \vec{d}^{n},p^{n+1}\right)  \in V_{E}\times Q_{D}$ such
that
\begin{align}
2\int_{D}\vec{d}^{n}\vec{v}+k_{n+1}\int_{D}\nu\,\nabla\vec{d}^{n} &  :\nabla
v+k_{n+1}\int_{D}\left(  \vec{w}^{n+1}\cdot\nabla\vec{d}^{n}\right)
v-\int_{D}p^{n+1}\left(  \nabla\cdot\vec{v}\right)  \label{eq:NS4-a}\\
&  =\int_{D}\frac{\partial\vec{u}^{n}}{\partial t}\vec{v}-\int_{D}\nu
\,\nabla\vec{u}^{n}:\nabla\vec{v}-\int_{D}\left(  \vec{w}^{n+1}\cdot\nabla
\vec{u}^{n}\right)  \vec{v},\nonumber\\
\int_{D}q\left(  \nabla\cdot\vec{d}^{n}\right)   &  =0,\label{eq:NS4-b}%
\end{align}
for all $(\vec{v},q)\in V_{D}\times Q_{D}$, and the TR velocity and
acceleration are updated as
\begin{equation}
\vec{u}^{n+1}=\vec{u}^{n}+k_{n+1}\vec{d}^{n},\qquad\frac{\partial\vec{u}%
^{n+1}}{\partial t}=2\vec{d}^{n}-\frac{\partial\vec{u}^{n}}{\partial
t}.\label{eq:NS4-c}%
\end{equation}

\paragraph{Time-step selection}

The time step size is driven by the heuristic formula
\begin{equation}
k_{n+2}=k_{n+1}\left(  \varepsilon/\left\Vert \vec{e}^{n+1}\right\Vert
\right)  ^{1/3}.\label{eq:time-step}%
\end{equation}
The local truncation error $\vec{e}^{n+1}$ is estimated by
\begin{equation}
\vec{e}^{n+1}=\left(  \vec{u}^{n+1}-\vec{u}_{\ast}^{n+1}\right)  /\left[
3\left(  1+k_{n}/k_{n+1}\right)  \right]  ,\label{eq:local-truncation-error}%
\end{equation}
where the TR\ velocity~$\vec{u}^{n+1}$ is compared with the AB2 velocity~$\vec
{u}_{\ast}^{n+1}$, which is computed using the explicit formula
\begin{equation}
\vec{u}_{\ast}^{n+1}=\vec{u}^{n}+\frac{k_{n+1}}{2}\left[  \left(
2+\frac{k_{n+1}}{k_{n}}\right)  \frac{\partial\vec{u}^{n}}{\partial t}-\left(
\frac{k_{n+1}}{k_{n}}\right)  \frac{\partial\vec{u}^{n-1}}{\partial t}\right]
.\label{eq:u-star}%
\end{equation}
There are three issues that need to be addressed:

\begin{enumerate}
\item \textit{The }AB2\textit{ is not self-starting}. To start the simulation
we require a function$~\vec{u}^{0}$ with boundary data$~\vec{g}^{0}$ such
that
\[
\int_{D}q\left(  \nabla\cdot\vec{u}^{0}\right)  =0,\quad\forall q\in Q_{D}.
\]
The initial acceleration (and pressure) is computed as follows: given the
boundary update $\vec{g}\vcentcolon=(\vec{g}^{1}-\vec{g}^{0})/k_{1}$, find the
pair $\left(  \frac{\partial\vec{u}^{0}}{\partial t},p^{0}\right)  \in
V_{E}\times Q_{D}$ such that
\begin{align*}
\int_{D}\frac{\partial\vec{u}^{0}}{\partial t}\vec{v}-\int_{D}p^{0}\left(
\nabla\cdot\vec{v}\right)   &  =-\int_{D}\nu\,\nabla\vec{u}^{0}:\nabla\vec
{v}-\int_{D}\left(  \vec{u}^{0}\cdot\nabla\vec{u}^{0}\right)  \vec{v},\\
\int_{D}q\left(  \nabla\cdot\frac{\partial\vec{u}^{0}}{\partial t}\right)   &
=0,
\end{align*}
for all $(\vec{v},q)\in V_{D}\times Q_{D}$. The discrete Oseen
problem~(\ref{eq:NS4-a})--(\ref{eq:NS4-b}) is then constructed by setting
$n=0$ and defining $\vec{w}^{1}=\vec{u}^{0}+k_{1}\frac{\partial\vec{u}^{0}%
}{\partial t}$, and its solution $\left(  \vec{u}^{1},p^{1}\right)  $ is used
to compute the acceleration at time $t=k_{1}$ as
\begin{equation}
\frac{\partial\vec{u}^{1}}{\partial t}=\frac{2}{k_{1}}\left(  \vec{u}^{1}%
-\vec{u}^{0}\right)  -\frac{\partial\vec{u}^{0}}{\partial t},
\label{eq:Oseen-acc-1}%
\end{equation}
and allows to compute the AB2 velocity at the second time step. The start-up
is completed by switching on the time-step control at the third time step
($k_{1}=k_{0}$).

\item \textit{Choice of initial time step}. The strategy is to select a
conservatively small value for$~k_{0}$, say $10^{-8}$. The time step then
typically exhibits a rapid growth in the first few steps, roughly as
$k_{n+1}/k_{n}=O\left(  \left(  \varepsilon/\mathtt{eps}\right)
^{1/3}\right)  \approx10^{4}$, with $\varepsilon=10^{-4}$ and considering the
(double) machine precision $\mathtt{eps}\approx10^{-16}$.

\item \textit{Time-step rejection}. The new time step is proposed by
formula~(\ref{eq:time-step}). However,\ if the next time step is seriously
reduced, i.e., $k_{n+2}<0.7k_{n+1}$ (or equivalently $\left\Vert \vec{e}%
^{n+1}\right\Vert >\left(  1/0.7\right)  ^{3}\varepsilon$), the next time step
is rejected: the value of $k_{n+1}$ is multiplied by $\left(  \varepsilon
/\left\Vert \vec{e}^{n+1}\right\Vert \right)  ^{1/3}$, and the current step is
repeated with this new$~k_{n+1}$.
\end{enumerate}

\paragraph{Stabilization of the integrator}

The numerical stabilization is implemented using \emph{time-step averaging}
with the purpose to annihilate any contribution of the form $(-1)^{n}$ to the
solution and its time derivative, which is invoked periodically every
$n_{\ast}$\ steps. For such a step the values of $t_{\ast}=t_{n}$ and $\vec
{u}^{\ast}=\vec{u}^{n}$\ are saved, we set $t_{n}=t_{n-1}+\frac{1}{2}k_{n}$,
$t_{n+1}=t_{\ast}+\frac{1}{2}k_{n+1}$ and define the new \textquotedblleft
shifted" solution vectors as
\begin{align*}
\vec{u}^{n} &  =\frac{1}{2}\left(  \vec{u}^{\ast}+\vec{u}^{n-1}\right)
,\qquad\frac{\partial\vec{u}^{n}}{\partial t}=\frac{1}{2}\left(
\frac{\partial\vec{u}^{n}}{\partial t}+\frac{\partial\vec{u}^{n-1}}{\partial
t}\right)  ,\\
\vec{u}^{n+1} &  =\vec{u}^{\ast}+\frac{1}{2}k_{n+1}\vec{d}^{n},\qquad
\frac{\partial\vec{u}^{n+1}}{\partial t}=\vec{d}^{n},
\end{align*}
where $\vec{d}^{n}$ is the TR\ update computed via~(\ref{eq:NS4-a}%
)--(\ref{eq:NS4-b}). In our implementation, the parameter$~n_{\ast}$ is fixed
and the value is set to$~10$.

\paragraph{Finite element formulation}

We consider the discretization of the Oseen problem~(\ref{eq:NS4-a}%
)--(\ref{eq:NS4-b}) by a div-stable mixed finite element method; in the
numerical experiments we use Taylor--Hood elements see,
e.g.,~\cite{Elman-2014-FEF}. Let the bases for velocity and pressure spaces be
denoted by $\left\{  \phi_{i}\right\}  _{i=1}^{n_{u}}$ and $\left\{
\varphi_{j}\right\}  _{j=1}^{n_{p}}$, respectively. In matrix terminology, the
Oseen problem at time step$~n$ entails solving a linear system
\begin{equation}
\left[
\begin{array}
[c]{cc}%
\mathbf{F}^{n+1} & \mathbf{B}^{T}\\
\mathbf{B} & \mathbf{0}%
\end{array}
\right]  \left[
\begin{array}
[c]{c}%
\mathbf{d}^{n}\\
\mathbf{p}^{n+1}%
\end{array}
\right]  =\left[
\begin{array}
[c]{c}%
\mathbf{f}_{v}^{n+1}\\
\mathbf{f}_{p}^{n+1}%
\end{array}
\right]  ,\label{eq:Oseen}%
\end{equation}
where $\mathbf{F}^{n+1}$ is the velocity convection-diffusion matrix: a sum of
the velocity mass matrix $\mathbf{M}$, diffusion matrix $\mathbf{A}$ and
convection matrix $\mathbf{N}^{n+1}$,\ defined as
\begin{equation}
\mathbf{F}^{n+1}=2\mathbf{M}+k_{n+1}\mathbf{A}+k_{n+1}\mathbf{N}%
^{n+1},\label{eq:F}%
\end{equation}
where
\begin{align*}
\mathbf{A}\mathbf{=}\left[  a_{ab}\right]  ,\qquad &  a_{ab}=\int_{D}%
\nu\,\nabla\phi_{b}:\nabla\phi_{a}\mathbf{,}\\
\mathbf{M}=\left[  m_{ab}\right]  ,\qquad &  m_{ab}=\int_{D}\phi_{b}\phi
_{a},\\
\mathbf{N}^{n+1}=\left[  n_{ab}^{n+1}\right]  ,\qquad &  n_{ab}^{n+1}=\int%
_{D}\left(  \vec{w}^{n+1}\cdot\nabla\phi_{b}\right)  \cdot\phi_{a},
\end{align*}
and $\vec{w}^{n+1}$ is computed from (\ref{eq:w}). The divergence
matrix~$\mathbf{B}$ is defined as
\begin{equation}
\mathbf{B}=\left[  b_{cd}\right]  ,\qquad b_{cd}=-\int_{D}\varphi_{c}\left(
\nabla\cdot\phi_{d}\right)  .\label{eq:B}%
\end{equation}
The right-hand side in~(\ref{eq:Oseen}) is constructed from the boundary
data$~\vec{g}^{n+1}$, the computed velocity$~\vec{u}^{n}$ at the previous time
level, and the acceleration$~\frac{\partial\vec{u}^{n}}{\partial t}$.

\section{Algorithms for the stochastic problem}

\label{sec:stochastic} Let $\left(  \Omega,\mathcal{F},\mathcal{P}\right)  $
represent a complete probability space, where $\Omega$ is the sample space,
$\mathcal{F}$ is a$~\sigma$-algebra on~$\Omega$ and $\mathcal{P}$ is a
probability measure. We assume that the randomness in the model is induced by
a vector $\xi:\Omega\rightarrow\Gamma\subset\mathbb{R}^{m_{\xi}}$ of
independent, identically distributed (i.i.d.)\ random variables $\xi
_{1}(\omega),\dots,\xi_{m_{\xi}}(\omega)$, where $\omega\in\Omega$. Let
$\mathcal{B}(\Gamma)$\ denote the Borel $\sigma$-algebra on $\Gamma$ induced
by$~\xi$, and $\mu$ denote the induced measure. The expected value of the
product of measurable functions on $\Gamma$\ determines a Hilbert
space~$T_{\Gamma}\equiv L^{2}\left(  \Gamma,\mathcal{B}(\Gamma),\mu\right)  $
with inner product
\begin{equation}
\left\langle u,v\right\rangle =\mathbb{E}\left[  uv\right]  =\int_{\Gamma
}u\left(  \xi\right)  v\left(  \xi\right)  \,d\mu\left(  \xi\right)  ,
\label{eq:E}%
\end{equation}
where the symbol~$\mathbb{E}$ denotes mathematical expectation.

In computations, we will use a finite-dimensional subspace $T_{p}\subset
T_{\Gamma}$ spanned by a set of multivariate polynomials $\left\{  \psi_{\ell
}(\xi)\right\}  $ that are orthonormal with respect to the density function
$\mu$, that is $\mathbb{E}\left[  \psi_{k}\psi_{\ell}\right]  =\delta_{k\ell}%
$, and $\psi_{1}=1$. This will be referred to as the gPC
basis~\cite{Xiu-2002-WAP}. The dimension of the space$~T_{p}$, depends on the
polynomial degree. For polynomials of total degree$~p_{\xi}$, the dimension is
$n_{\xi}=\binom{m_{\xi}+p_{\xi}}{p_{\xi}}$.

\subsection{Navier\textendash Stokes equation with stochastic viscosity}

We use the same setup as in~\cite{Sousedik-2016-SGM}. Specifically, we
consider that the expansion of viscosity is given as
\begin{equation}
\nu\equiv\nu\left(  x,\xi\right)  =\sum_{\ell=1}^{n_{\nu}}\nu_{\ell}%
(x)\psi_{\ell}(\xi),\label{eq:KL}%
\end{equation}
where $\nu_{\ell}(x)$ is a set of deterministic spatial functions,  
\textcolor{black}{and index~$\ell$ is related through a multi-index to the degrees of the random 
variables $\xi_1,\dots,\xi_{m_\xi}$ used in the construction of the gPC basis function~$\psi_{\ell}(\xi)$ see,
e.g.,~\cite[Section~2.4.3]{Ghanem-1991-SFE} or~\cite[Section~5.2]{Xiu-2010-NMS}.}
For simplicity, we will also assume that both the Dirichlet boundary
conditions~(\ref{eq:bc-Dirichlet}) and the initial condition~(\ref{eq:ic}) are
deterministic. We seek a discrete approximation of the velocity in the form
\begin{equation}
\vec{u}\left(  x,t,\xi\right)  \approx\sum_{k=1}^{n_{\xi}}\sum_{i=1}^{n_{u}%
}u_{ik}(x,t)\phi_{i}(x)\psi_{k}(\xi)=\sum_{k=1}^{n_{\xi}}\vec{u}_{k}%
(x,t)\psi_{k}(\xi).\label{eq:u-gPC}%
\end{equation}

\begin{remark}
In literature it is sometime recommended to use a time-dependent gPC\ basis,
that is $\psi_{k}(\xi,t)$, to keep the stochastic dimension low in long-time
integration. However, this is a complementary strategy to the present study.
Since it is not needed in our numerical experiments, we use only a
time-independent gPC\ basis.
\end{remark}

\subsection{Stochastic Galerkin method}

\label{sec:SGM}The stochastic Galerkin formulation of problem (\ref{eq:NS3-a}%
)--(\ref{eq:NS3-b}) consists of using the expansion~(\ref{eq:KL}) and
performing a Galerkin projection on the space~$T_{\Gamma}$\ using mathematical
expectation in the sense of~(\ref{eq:E}). That is, we seek velocity $\vec
{u}^{n+1}\in T_{\Gamma}\otimes V_{E}$ and pressure $p^{n+1}\in T_{\Gamma
}\otimes Q_{D}$ for a given pair $(\vec{u}^{n},p^{n})$, such that
\begin{align*}
\mathbb{E}\left[  \frac{2}{k_{n+1}}\int_{D}\vec{u}^{n+1}\vec{v}+\int_{D}%
\nu\,\nabla\vec{u}^{n+1}\right.   &  \left.  :\nabla\vec{v}+\int_{D}\left(
\vec{w}^{n+1}\cdot\nabla\vec{u}^{n+1}\right)  \vec{v}-\int_{D}p^{n+1}\left(
\nabla\cdot\vec{v}\right)  \right] \\
&  =\mathbb{E}\left[  \frac{2}{k_{n+1}}\int_{D}\vec{u}^{n}\vec{v}+\int%
_{D}\frac{\partial\vec{u}^{n}}{\partial t}\vec{v}\right]  ,\quad\forall\vec
{v}\in T_{\Gamma}\otimes V_{D},\\
\mathbb{E}\left[  \int_{D}q\left(  \nabla\cdot\vec{u}^{n+1}\right)  \right]
&  =0,\quad\forall q\in T_{\Gamma}\otimes Q_{D},
\end{align*}
and the stochastic counterpart of the discrete Oseen problem~(\ref{eq:NS4-a}%
)--(\ref{eq:NS4-c}) is: given $\vec{u}^{n}$, $\partial\vec{u}^{n}/\partial
t$\ and the boundary update $\vec{g}\vcentcolon=(\vec{g}^{n+1}-\vec{g}%
^{n})/k_{n+1}$, we first compute $\vec{d}^{n}\in T_{\Gamma}\otimes V_{E}$ and
$p^{n+1}\in T_{\Gamma}\otimes Q_{D}$\ such that
\begin{align}
\mathbb{E}\left[  2\int_{D}\vec{d}^{n}\vec{v}\right.   &  +\left.  k_{n+1}%
\int_{D}\nu\,\nabla\vec{d}^{n}:\nabla v+k_{n+1}\int_{D}\left(  \vec{w}%
^{n+1}\cdot\nabla\vec{d}^{n}\right)  v-\int_{D}p^{n+1}\left(  \nabla\cdot
\vec{v}\right)  \right] \label{eq:Oseen-stoch-1}\\
&  =\mathbb{E}\left[  \int_{D}\frac{\partial\vec{u}^{n}}{\partial t}\vec
{v}-\int_{D}\nu\,\nabla\vec{u}^{n}:\nabla\vec{v}-\int_{D}\left(  \vec{w}%
^{n+1}\cdot\nabla\vec{u}^{n}\right)  \vec{v}\right]  ,\quad\forall\vec{v}\in
T_{\Gamma}\otimes V_{D},\nonumber\\
\mathbb{E}\left[  \int_{D}q\left(  \nabla\cdot\vec{d}^{n}\right)  \right]   &
=0,\quad\forall q\in T_{\Gamma}\otimes Q_{D}, \label{eq:Oseen-stoch-2}%
\end{align}
and the TR velocity and the acceleration are updated as in~(\ref{eq:NS4-c}).

\subsubsection{Stochastic Galerkin finite element formulation}

The Galerkin projection leads to a large coupled system of equations with
structure depending on the ordering of the unknown coefficientsts~$\{u_{ik}%
\}$, $\{p_{jk}\}$. We will group velocity-pressure pairs for each~$k$, the
index of stochastic basis functions (and order equations in the same way),
giving the ordered list of coefficients
\begin{equation}
u_{1:n_{u},1},p_{1:n_{p},1},\,u_{1:n_{u},2},p_{1:n_{p},2},\,\ldots
,\,u_{1:n_{u},n_{\xi}},p_{1:n_{p},n_{\xi}}.\label{eq:coefficient-order}%
\end{equation}
The discrete stochastic Oseen operator is built as follows. First, we set up
the discrete components of the diffusion matrix using the expansion of
viscosity~(\ref{eq:KL}) as
\begin{equation}
\mathbf{A}_{\ell}\mathbf{=}\left[  a_{\ell,ab}\right]  ,\quad a_{\ell
,ab}=\left(  \int_{D}\nu_{\ell}(x)\,\nabla\phi_{b}:\nabla\phi_{a}\right)
,\qquad\ell=1,\ldots,n_{\nu}.\label{stochastic-vector-Laplacian}%
\end{equation}
Next, let~$\vec{w}_{\ell}^{n+1}(x)$ denote the $\ell$th term of the
extrapolated velocity iterate (as in the expression on the right
in~(\ref{eq:u-gPC}) for $k=\ell$) at step~$n$, and let
\[
\mathbf{N}_{\ell}^{n+1}=\left[  n_{\ell,ab}^{n+1}\right]  ,\qquad n_{\ell
,ab}^{n+1}=\int_{D}\left(  \vec{w}_{\ell}^{n+1}\cdot\nabla\phi_{b}\right)
\cdot\phi_{a},\qquad\ell=1,\ldots,n_{\xi}.
\]
Let $\widehat{n}=\max(n_{\nu},n_{\xi})$ and, if needed, define $\mathbf{A}%
_{\ell}=\mathbf{0}$ for $n_{\nu}<\ell\leq\widehat{n}$\ and $\mathbf{N}_{\ell
}^{n+1}=\mathbf{0}$\ for $n_{\xi}<\ell\leq\widehat{n}$. Then in analogue
to~(\ref{eq:F}) define matrices
\begin{align}
\mathbf{F}_{1}^{n+1} &  =2\mathbf{M}+k_{n+1}\mathbf{A}_{1}+k_{n+1}%
\mathbf{N}_{1}^{n+1},\qquad\label{eq:F_1-s}\\
\mathbf{F}_{\ell}^{n+1} &  =k_{n+1}\mathbf{A}_{\ell}+k_{n+1}\mathbf{N}_{\ell
}^{n+1},\quad\ell=2,\dots,\widehat{n},\label{eq:F_2-s}%
\end{align}
which are incorporated into the block matrices
\begin{equation}
\mathcal{F}_{1}^{n+1}=\left[
\begin{array}
[c]{cc}%
\mathbf{F}_{1}^{n+1} & \mathbf{B}^{T}\\
\mathbf{B} & \mathbf{0}%
\end{array}
\right]  ,\qquad\mathcal{F}_{\ell}^{n+1}=\left[
\begin{array}
[c]{cc}%
\mathbf{F}_{\ell}^{n+1} & \mathbf{0}\\
\mathbf{0} & \mathbf{0}%
\end{array}
\right]  ,\quad\ell=2,\dots,\widehat{n}.\label{eq:stochastic-saddle-point}%
\end{equation}
These operators will be coupled with matrices arising from terms in~$T_{p}$,
\begin{equation}
\mathbf{H}_{\ell}=\left[  h_{\ell,jk}\right]  ,\quad h_{\ell,jk}%
\equiv\mathbb{E}\left[  \psi_{\ell}\psi_{j}\psi_{k}\right]  ,\qquad
\ell=1,\dots,n_{\nu},\quad j,k=1,\dots,n_{\xi}%
.\label{eq:stochastic-matrix-forms}%
\end{equation}
Combining the expressions from~(\ref{eq:stochastic-saddle-point})
and~(\ref{eq:stochastic-matrix-forms}), using the
ordering~(\ref{eq:coefficient-order}) yields the discrete stochastic Oseen
system
\begin{equation}
\left(  \sum_{\ell=1}^{\widehat{n}}\mathbf{H}_{\ell}\otimes\mathcal{F}_{\ell
}^{n+1}\right)  \mathbf{v}=\mathbf{y},\label{eq:Oseen-s}%
\end{equation}
where $\otimes$ denotes the matrix Kronecker product. The entries of the
vectors$~\mathbf{v}$ and$\mathbf{~y}$ are ordered as
in~(\ref{eq:coefficient-order}). Note that~$\mathbf{H}_{1}$ is the identity
matrix of order~$n_{\xi}$.

\begin{remark}
\label{rem:order} With this ordering, which we used also
in~\cite{Sousedik-2016-SGM}, the coefficient matrix contains a set of~$n_{\xi
}$ block~$2\times2$ matrices of saddle-point structure along its block
diagonal, given by
\[
\mathcal{F}_{1}^{n+1}+\sum_{\ell=2}^{\widehat{n}}h_{\ell,jj}\mathcal{F}_{\ell
}^{n+1},\qquad j=1,\ldots,n_{\xi}.
\]
This enables the use of existing deterministic solvers for the individual
diagonal blocks.

\end{remark}

We find it convenient to formulate the solvers in the so-called matricized
format. To this end, we make use of isomorphism between~$%
\mathbb{R}
^{n_{x}n_{\xi}}$ and $%
\mathbb{R}
^{n_{x}\times n_{\xi}}$\ determined by the operators $\operatorname{vec}$ and
$\operatorname{mat}$. Let $n_{x}=n_{u}+n_{p}$ and consider writing the
solution of~(\ref{eq:Oseen-s}) using the ordering~(\ref{eq:coefficient-order})
as $\mathbf{v}=[v_{1}^{T},v_{2}^{T},\dots,v_{n_{\xi}}^{T}]^{T}$, where
$v_{k}=[\vec{u}_{k}^{T},p_{k}^{T}]$\ for $k=1,\dots,n_{\xi}$\ as in the
expansions on the right in (\ref{eq:u-gPC}).
Then we write $\mathbf{v}=\operatorname{vec}(\mathbf{V})$, $\mathbf{V}%
=\operatorname{mat}(\mathbf{v})$, where $\mathbf{v}\in%
\mathbb{R}
^{n_{x}n_{\xi}}$, $\mathbf{V}\in%
\mathbb{R}
^{n_{x}\times n_{\xi}}$ and the upper/lower case notation is assumed
throughout the paper, so $\mathbf{Y}=\operatorname{mat}(\mathbf{y})$, etc.
Specifically, we define the \emph{matricized}\ coefficients of the solution
expansion
\begin{equation}
\mathbf{V}=\operatorname{mat}(\mathbf{v})=\left[  v_{1},v_{2},\ldots
,v_{n_{\xi}}\right]  \in%
\mathbb{R}
^{n_{x}\times n_{\xi}},\label{eq:U}%
\end{equation}
where the column~$k$ contains the coefficients associated with the basis
function$~\psi_{k}$. In this setting, since $\left(  \mathbf{V}\otimes
\mathbf{W}\right)  \operatorname{vec}\left(  \mathbf{X}\right)
=\operatorname{vec}\left(  \mathbf{WXV}^{T}\right)  $, the linear
system~(\ref{eq:Oseen-s}) can be equivalently written as
\begin{equation}
\sum_{\ell=1}^{\widehat{n}}\mathcal{F}_{\ell}^{n+1}\mathbf{VH}_{\ell
}=\mathbf{Y}.\label{eq:Oseen-s-mat}%
\end{equation}

The time-step selection is driven by the formula~(\ref{eq:time-step}), which
we heuristically modify as follows. First, we run the deterministic solver
with viscosity $\nu=\nu_{1}$ and record the set of time steps $0,t_{1}%
,t_{2},,\dots,T$. Then, we divide size of each interval $[t_{n},t_{n+1}]$
by~$n_{\xi}$, and we further round down the time-step size to the nearest
power of~$10$. This procedure yields a sequence of time steps, which is then
used for evolution of the stochastic Galerkin method.
We note that an alternative strategy could be utilized by using the
gPC\ coefficients corresponding to the mean velocity directly in
formula~(\ref{eq:time-step}), that is without an a priori run of the
deterministic solver.

\subsection{Sampling methods}

\label{sec:sampling}
Both Monte Carlo and stochastic collocation methods are based on sampling.
This entails the solution of a number of mutually independent deterministic
problems at a set of sample points~$\left\{  \xi^{\left(  q\right)  }\right\}
$, which give realizations of the viscosity~(\ref{eq:KL}). That is, a
realization of viscosity~$\nu\left(  \xi^{\left(  q\right)  }\right)  $ gives
rise to deterministic functions~$\vec{u}\left(  \cdot,\cdot,\xi^{\left(
q\right)  }\right)  $ and$~p\left(  \cdot,\cdot,\xi^{\left(  q\right)
}\right)  $ on~$D$ that satisfy the standard deterministic Navier\textendash
Stokes equations, and to corresponding finite-element approximations.

In the Monte Carlo method, the$~n_{MC}$ sample points are generated randomly,
following the distribution of the random variables$~\xi$, and moments of the
solution are obtained from ensemble averaging. In addition the coefficients
in~(\ref{eq:u-gPC}) could be determined at time~$t_{b}$ using\footnote{\textcolor{black}{In numerical experiments, 
we avoid this approximation of the gPC coefficients and directly work with the sampled quantities.}}
\[
u_{ik}(t_{b})=\frac{1}{n_{MC}}\sum_{q=1}^{n_{MC}}\vec{u}^{\,(q)}\left(
x_{i},t_{b}\right)  \,\psi_{k}\left(  \xi^{\left(  q\right)  }\right)  ,
\]
where for$~t_{b}$, $b=1,\dots,n_{b}$, we will consider an a priori set of
\emph{time barriers}, which is used in implementation to enforce 
all~$n_{MC}$\ instances of the deterministic solver to step through. For
stochastic collocation, the sample points consist of a set of predetermined
\emph{collocation points}. This approach derives from a methodology for
performing quadrature or interpolation in multidimensional space using a small
number of points, a so-called sparse
grid~\cite{Gerstner-1998-NIU,Novak-1996-HDI}. There are \textcolor{black}{several} ways to implement
stochastic collocation to obtain the coefficients in~(\ref{eq:u-gPC}). 
\textcolor{black}{In the basic variant of the method, it is possible proceed}
either by constructing a Lagrange interpolating polynomial, or, in the
so-called pseudospectral approach, by performing a discrete
projection~\cite{Xiu-2010-NMS}. 
We use the \textcolor{black}{pseudospectral} 
approach because it
facilitates a direct comparison with the stochastic Galerkin method, and we
refer, e.g., to~\cite{LeMaitre-2010-SMU} for an overview and discussion of
integration rules. In particular, the coefficients in~(\ref{eq:u-gPC}) are
determined at time~$t_{b}$ using a quadrature rule
\[
u_{ik}(t_{b})=\sum_{q=1}^{n_{q}}\vec{u}^{\,(q)}\left(  x_{i},t_{b}\right)
\,\psi_{k}\left(  \xi^{\left(  q\right)  }\right)  \,w^{\left(  q\right)  },
\]
where~$\xi^{\left(  q\right)  }$ and~$w^{\left(  q\right)  }$, $q=1,\dots
,n_{q}$, are the collocation (quadrature) points and weights.
\textcolor{black}{Finally, we note that the other ways to 
perform stochastic collocation include the least-square approach and the compressed sensing approach see,
e.g.,~\cite{Chkifa-2015-DLS,Guo-2020-CLP,Hampton-2015-CSP,Narayan-2017-CFW}.}

\section{Numerical experiments}

\label{sec:numerical}
\begin{figure}[ptbh]
\begin{center}
\includegraphics[width=12.8cm]{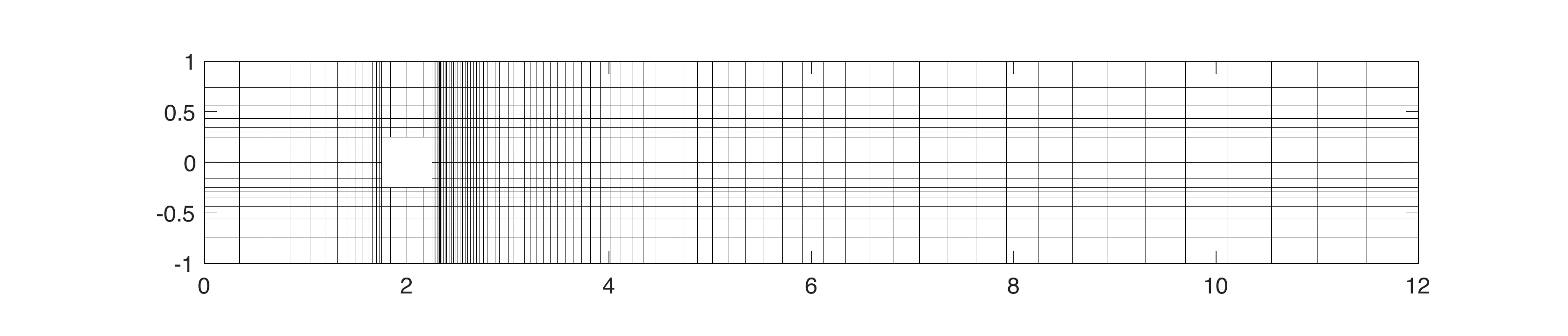} \vspace{-4mm}
\end{center}
\caption{Finite element mesh for the flow around an obstacle problem. }%
\label{fig:mesh-obstacle}%
\end{figure}

We implemented the method in \textsc{Matlab} using the
IFISS~3.5 package~\cite{Elman-2014-IFISS,ifiss}, and in this section we
present results of numerical experiments for a model problem given by a flow
around an obstacle. The geometry of the problem is shown in
Figure~\ref{fig:mesh-obstacle}. The discretization of the physical space
consists of $12,640$ velocity and $1640$ pressure degrees of freedom. The
viscosity was taken to be a lognormal process, and its representation was
computed from an underlying Gaussian random process using the transformation
described in \cite{Ghanem-1999-NGS}. That is, for $\ell=1,\dots,n_{\nu}$,
$\psi_{\ell}\left(  \xi\right)  $ is the product of$~m_{\xi}$ univariate
Hermite polynomials, and denoting the coefficients of the Karhunen-Lo\`{e}ve
expansion of the Gaussian process by~$g_{j}\left(  x\right)  $ and $\eta
_{j}=\xi_{j}-g_{j}$, $j=1,\dots,m_{\xi}$, the coefficients in
expansion~(\ref{eq:KL}) are computed as
\[
\nu_{\ell}\left(  x\right)  =\mathbb{E}\left[  \psi_{\ell}\left(  \eta\right)
\right]  \exp\left[  g_{0}\left(  x\right)  +\frac{1}{2}\sum_{j=1}^{m_{\xi}%
}\left(  g_{j}\left(  x\right)  \right)  ^{2}\right]  .
\]
The covariance function of the Gaussian field, for points~$X_{i}=(x_{i}%
,y_{i})\in D$, $i=1,2$, was chosen to be
\[
C\left(  X_{1},X_{2}\right)  =\sigma_{g}^{2}\exp\left(  -\frac{\left\vert
x_{2}-x_{1}\right\vert }{L_{x}}-\frac{\left\vert y_{2}-y_{1}\right\vert
}{L_{y}}\right)  ,
\]
where~$L_{x}$ and $L_{y}$\ are the correlation lengths of the random
variables~$\xi_{i}$, $i=1,\dots,m_{\xi}$, in the $x$ and $y$ directions,
respectively, and $\sigma_{g}$ is the standard deviation of the Gaussian
random field. The correlation lengths were set to be equal to $25\%$ of the
width and height of the domain, i.e. $L_{x}=3$ and $L_{y}=0.5$. The
coefficient of variation of the lognormal field, defined as $CoV=\sigma_{\nu
}/\nu_{1}$ where $\sigma_{\nu}$ is the standard deviation, was set to either
$1\%$ or $10\%$. The stochastic dimension was $m_{\xi}=2$. The degree used for
the polynomial expansion of the solution was $p_{\xi}=3$, and the degree used
for the expansion of the lognormal process was $2p_{\xi}$, which ensures a
complete representation of the process in the discrete
problem~\cite{Matthies-2005-GML}. With these settings, $n_{\xi}=10$ and
$n_{\nu}=\widehat{n}=28$, and $\mathbf{H}_{\ell}$ is of order $10$
in~(\ref{eq:Oseen-s}). For the mean value of viscosity we used $\nu_{1}=0.02$,
which corresponds to mean Reynolds number $\operatorname*{Re}_{1}=100$, and
$\nu_{1}=6.67\times10^{-3}$, which corresponds to mean Reynolds number
$\operatorname*{Re}_{1}=300$. We note that the steady-state case was studied
in~\cite{Sousedik-2016-SGM}, and the setup for the (deterministic)
time-dependent problem is the same to~\cite[Chapter~10]{Elman-2014-FEF} except
that the length of the channel was set to~$12$. Specifically, the initial
condition for velocity was taken zero, and the Dirichlet boundary condition on
the inflow$~\partial D_{\text{Dir}}$ (the left side) was smoothly ramped up
from zero to steady state as $\vec{u}\left(  \cdot,t\right)  =\left(
1-e^{-5t}\right)  \vec{w},$ where $\vec{w}$ is a Poiseuille (parabolic) flow
profile, no-flow condition was prescribed on the top and bottom walls and
natural `do-nothing' condition was used on the outflow boundary (the right
side). The initial time step was set to$~k_{0}$=$10^{-9}$, and the problem was
evolved from $t_{0}=0\,s$ to $T=10\,s$. The time-stepping method described in
Section~\ref{sec:model} was used for each sample of viscosity $\nu\left(
x,\xi^{(q)}\right)  $ for both Monte Carlo and the stochastic collocation
methods, and the method from Section~\ref{sec:SGM} was used for the stochastic
Galerkin method. In order to compare the gPC\ coefficients at the a priori
chosen set of times, which we will refer to as \emph{time barriers}, we
prescribed the stochastic Galerkin solver to step through certain times, and
we used the same set also for the Monte Carlo simulation in order to compare
probability density function estimates of the velocity obtained by using all
three methods. Specifically, we used time barriers $t_{b}%
=\{0,0.1,0.2,0.5,1,2,5,6,8,10\}$.

\begin{figure}[ptbh]
\begin{center}
\includegraphics[width=10cm]{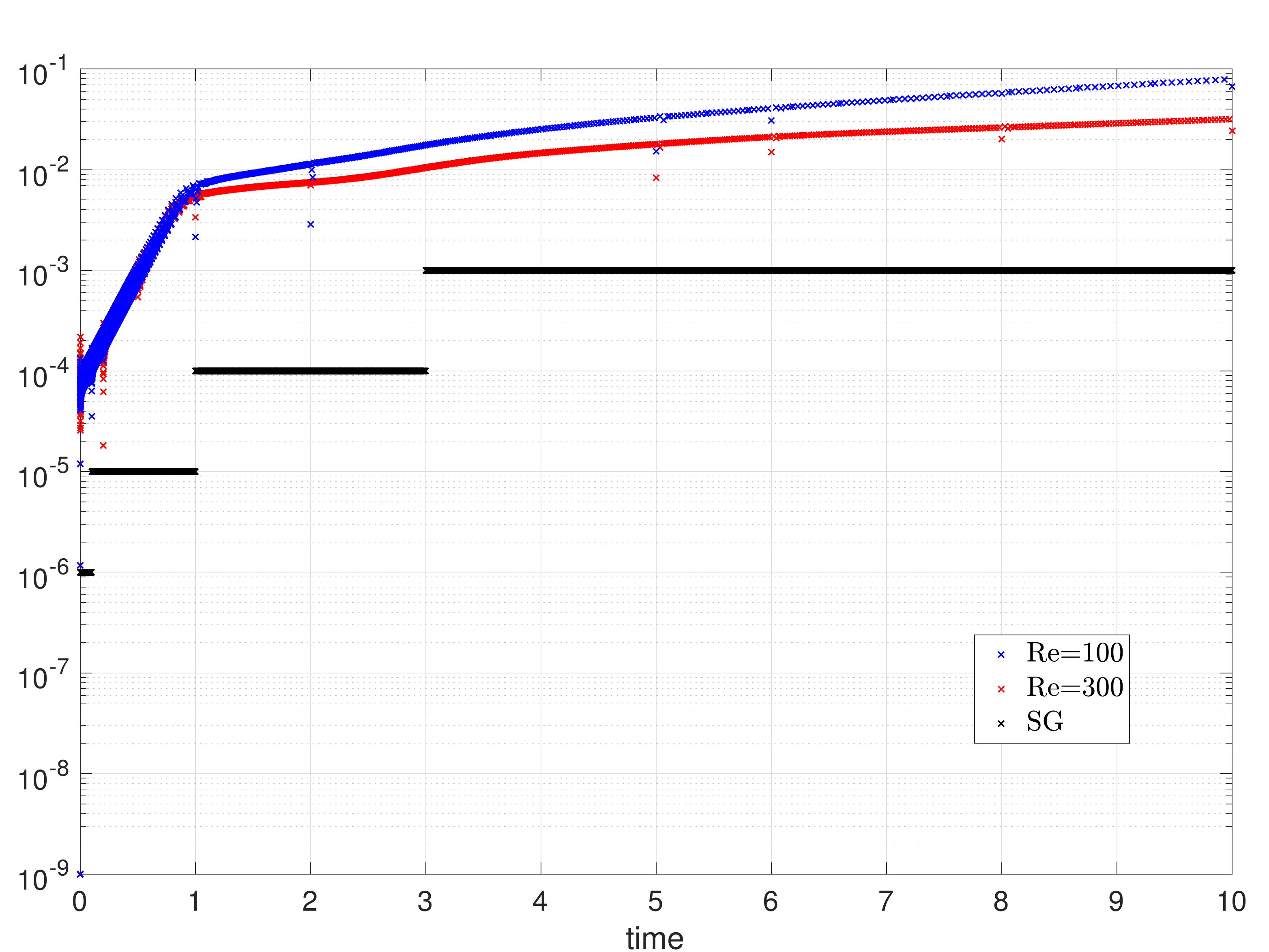}
\end{center}
\par
\caption{Evolution of the time-step size for the deterministic problems with
Reynolds numbers $\operatorname*{Re}_{1}=100$ and $\operatorname*{Re}_{1}=300$
and for the stochastic Galerkin method (SG).}%
\label{fig:timestepplot}%
\end{figure}

The evolution of the time step for the two deterministic cases with mean
Reynolds numbers $\operatorname*{Re}_{1}=100$ and $\operatorname*{Re}_{1}=300$
and for the stochastic Galerkin method is shown in
Figure~\ref{fig:timestepplot}. We note that the heuristic used for the
stochastic Galerkin methods yields the same time-step selection for both
values of the Reynolds number. Specifically, only three steps with size
$10^{-9}$ are performed at the very beginning, then the step size increases to
$10^{-6}$ and eventually to $10^{-5}$ for the most part of the first second.
For the next two seconds it becomes $10^{-4}$ and eventually $10^{-3}$ for the
rest of the time.

\begin{figure}[ptbh]
\begin{center}
\includegraphics[width=10cm]{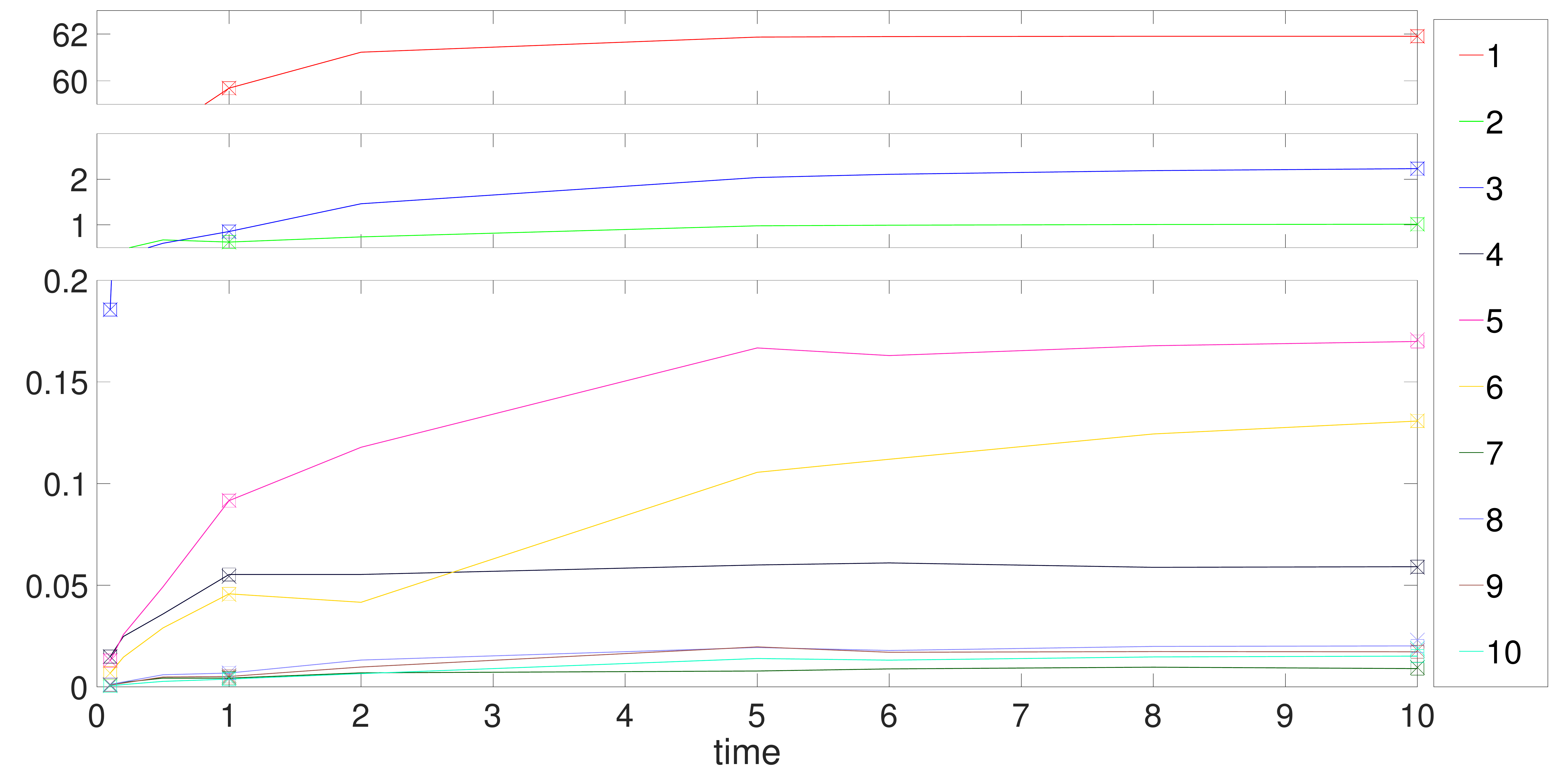}
\end{center}
\caption{Evolution of the gPC coefficients corresponding to the horizontal
velocity in terms of $\ell_{2}$-norm for mean Reynolds number
$\operatorname*{Re}_{1}=100$ and $CoV=10\%$. The symbols $\square$ and
$\times$ represent the results of the Monte Carlo and stochastic collocation,
respectively, at times $0.1$s, $1$s and $10$s.}%
\label{fig:gPC_100_10}%
\end{figure}

\begin{figure}[ptbh]
\begin{center}
\includegraphics[width=12cm]{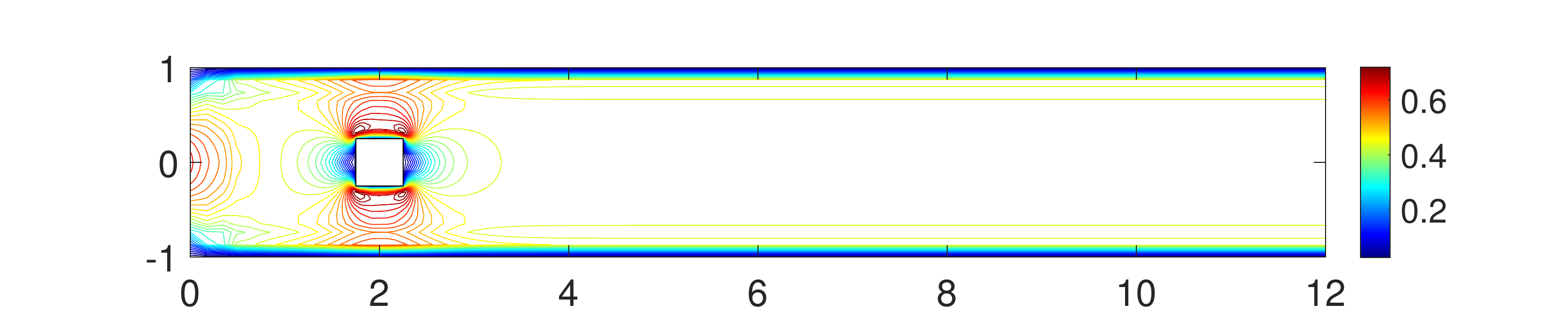}
\includegraphics[width=12cm]{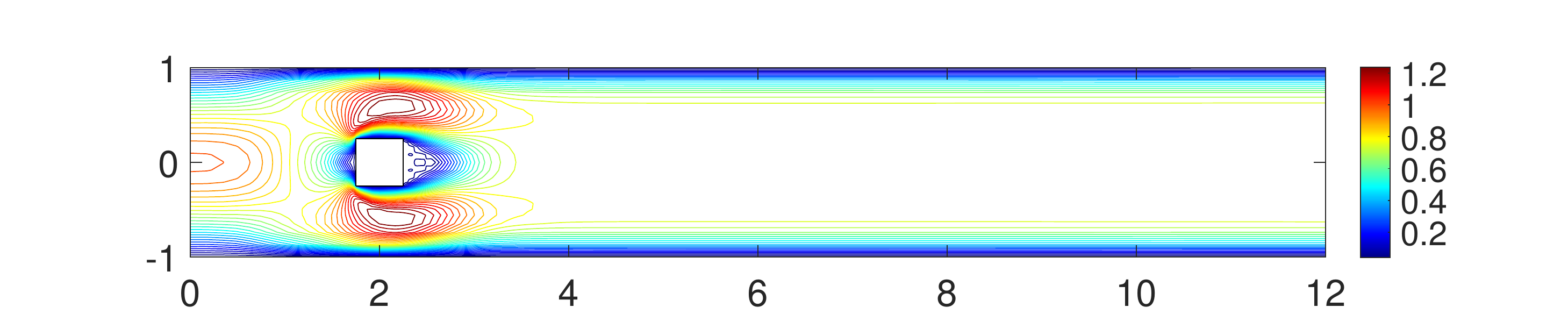}
\includegraphics[width=12cm]{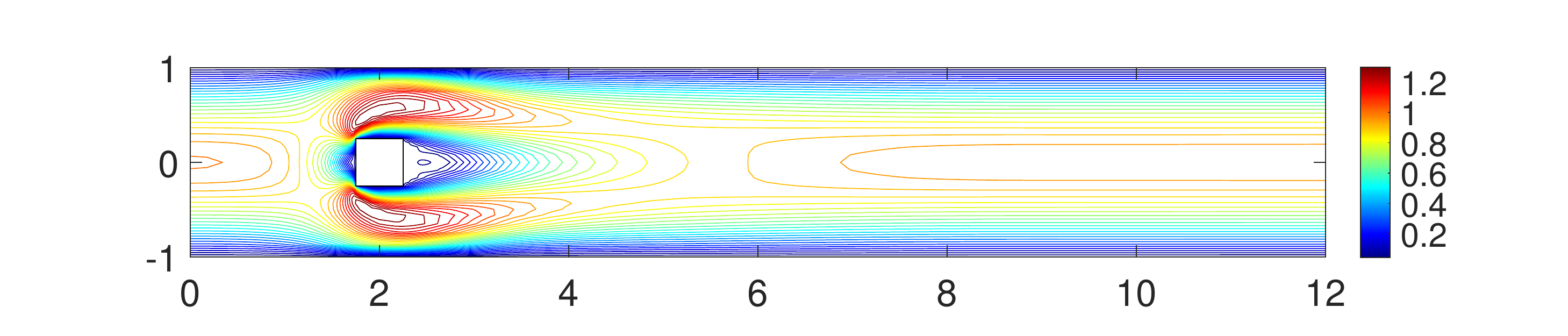}
\end{center}
\caption{Mean horizontal velocity at times $0.1$s (top), $1$s (center) and
$10$s (bottom) for mean Reynolds number $\operatorname*{Re}_{1}=100$ and
$CoV=10\%$.}%
\label{fig:mean_velocity_100_10}%
\end{figure}

\begin{figure}[ptbh]
\begin{center}
\includegraphics[width=12cm]{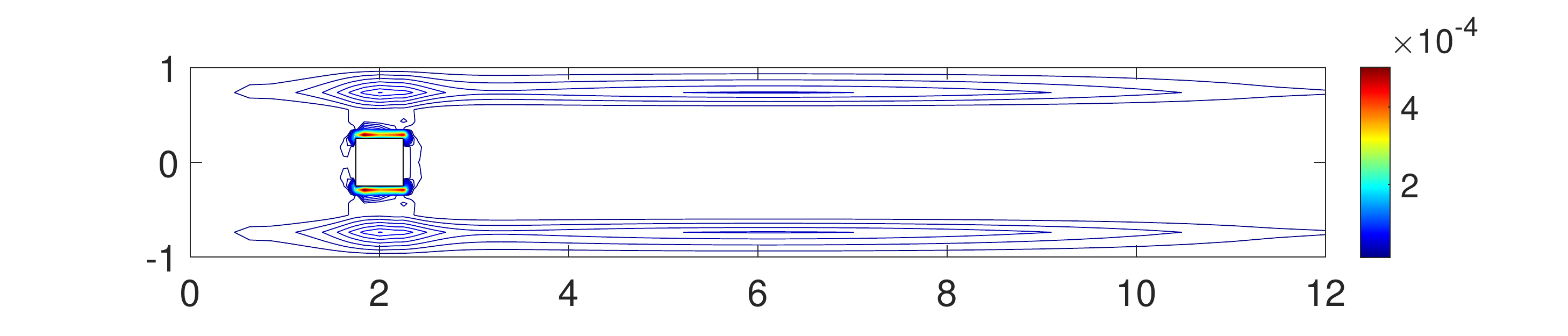}
\includegraphics[width=12cm]{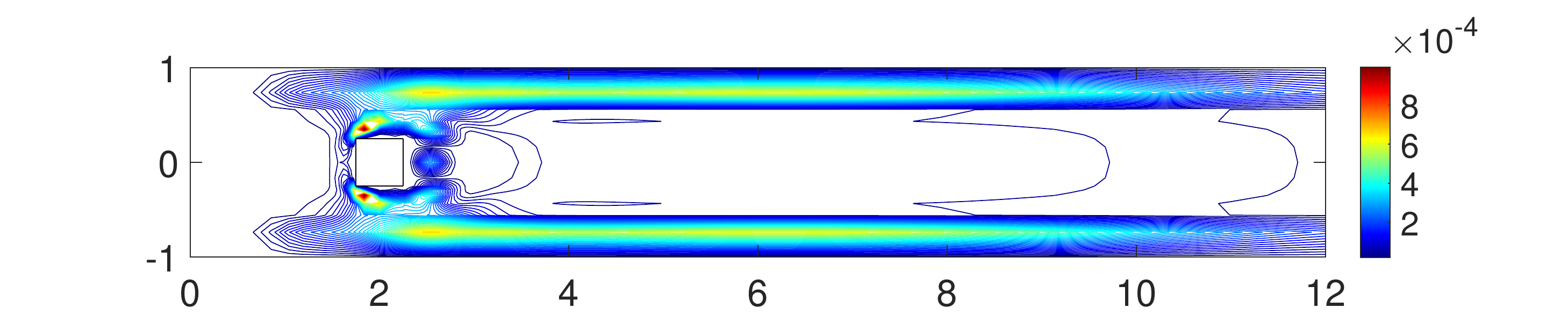}
\includegraphics[width=12cm]{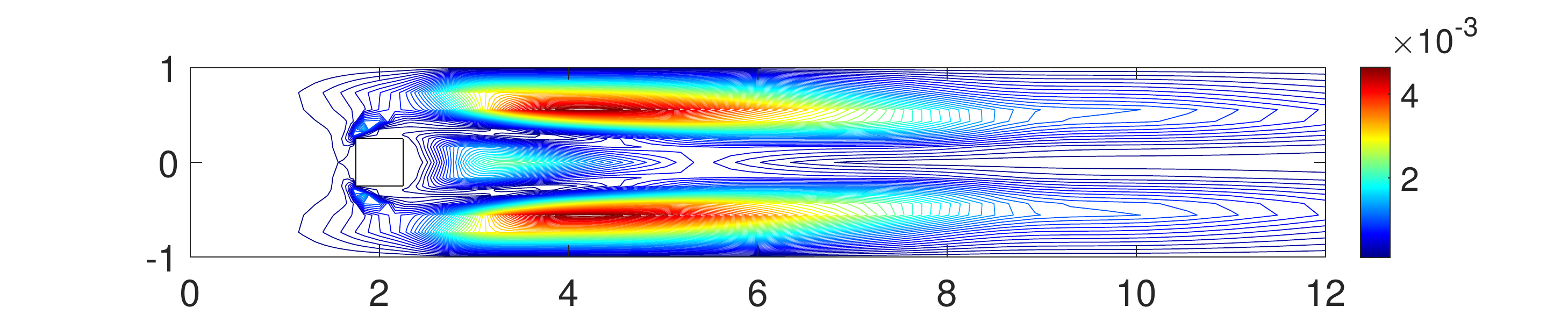}
\end{center}
\caption{Variance of the horizontal velocity at times $0.1$s (top), $1$s
(center) and $10$s (bottom) for mean Reynolds number $\operatorname*{Re}%
_{1}=100$ and $CoV=10\%$.}%
\label{fig:variance_ux_100_10}%
\end{figure}

\begin{figure}[ptbh]
\begin{center}
\includegraphics[width=12cm]{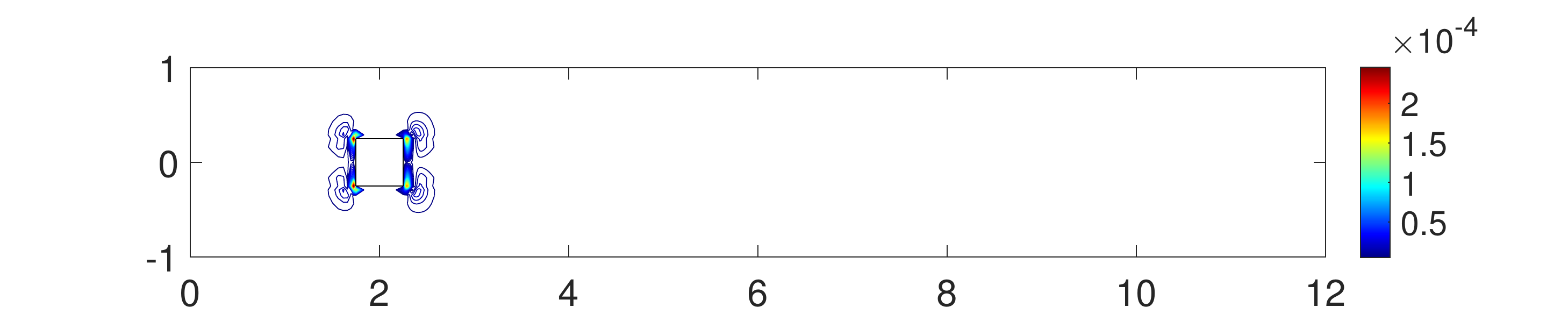}
\includegraphics[width=12cm]{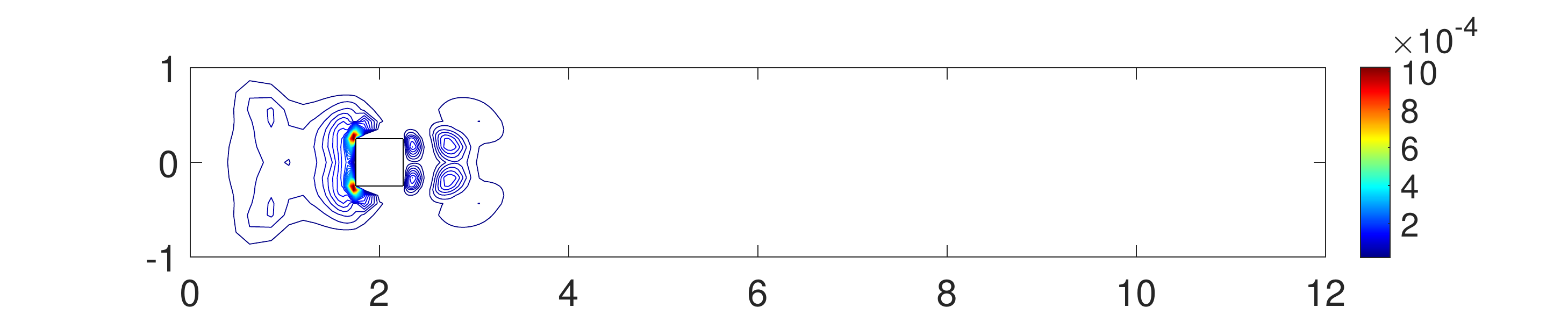}
\includegraphics[width=12cm]{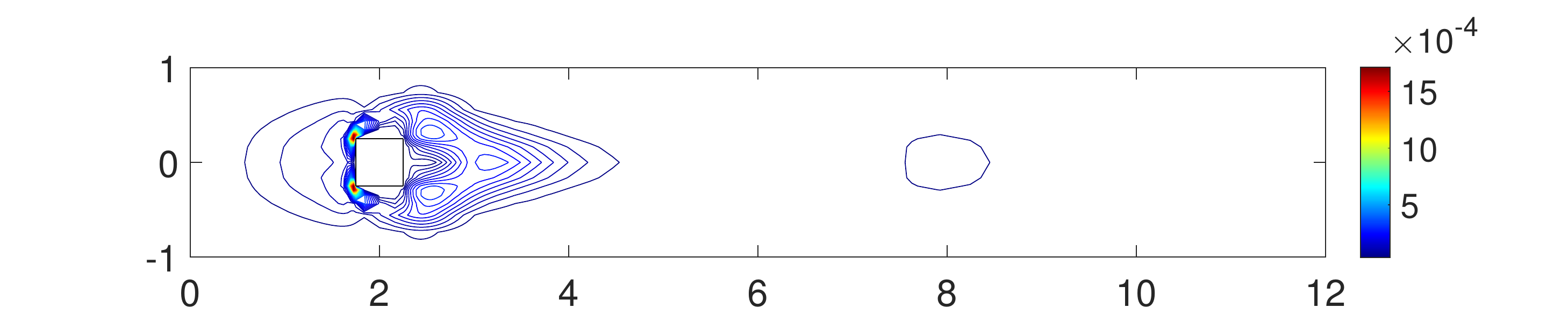}
\end{center}
\caption{Variance of the vertical velocity at times $0.1$s (top), $1$s
(center) and $10$s (bottom) for mean Reynolds number $\operatorname*{Re}%
_{1}=100$ and $CoV=10\%$.}%
\label{fig:variance_uy_100_10}%
\end{figure}

\begin{figure}[ptbh]
\begin{center}%
\begin{tabular}
[c]{c@{\hskip-0.2cm}c}%
\includegraphics[width=6.5cm]{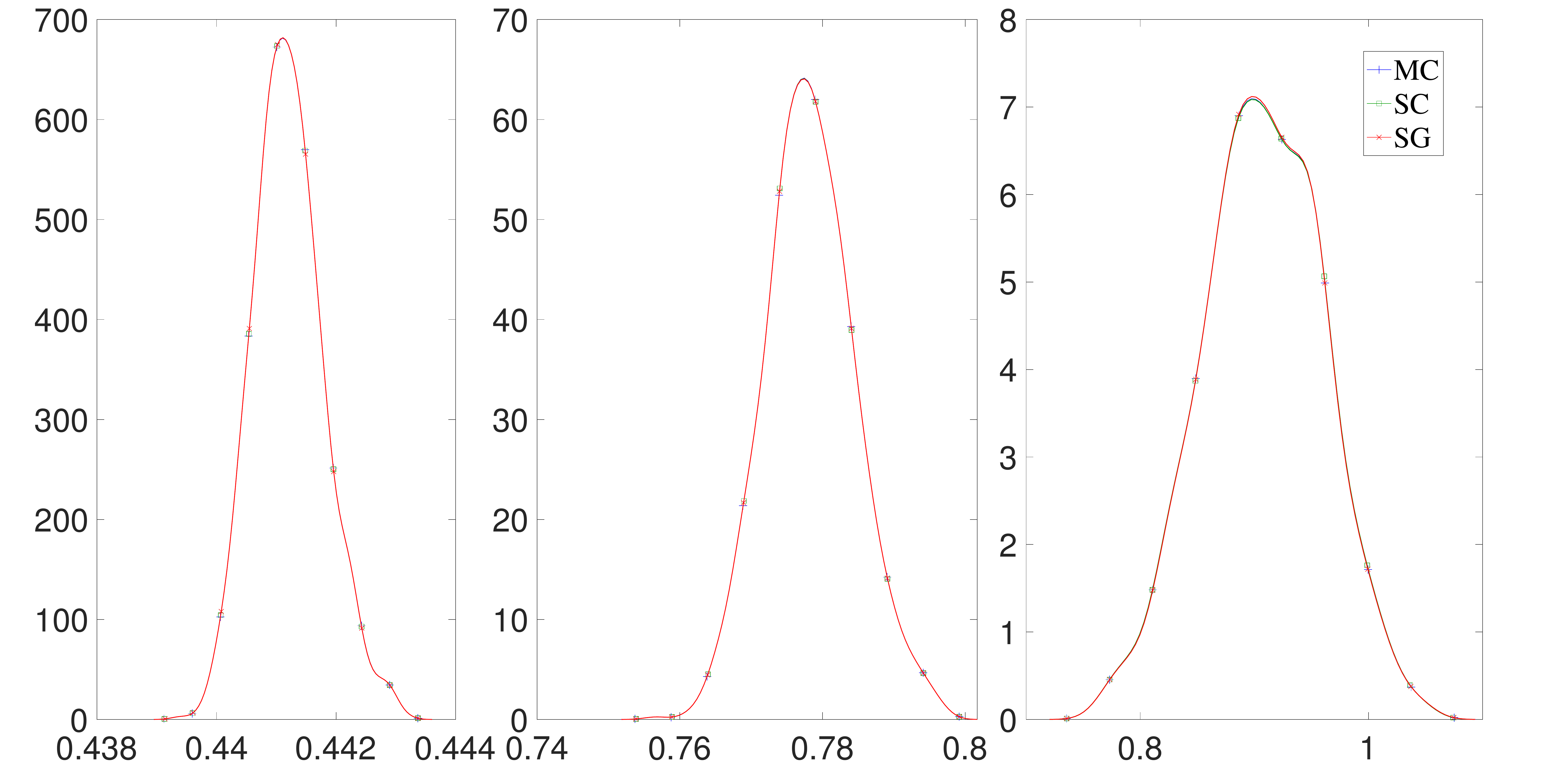} &
\includegraphics[width=6.5cm]{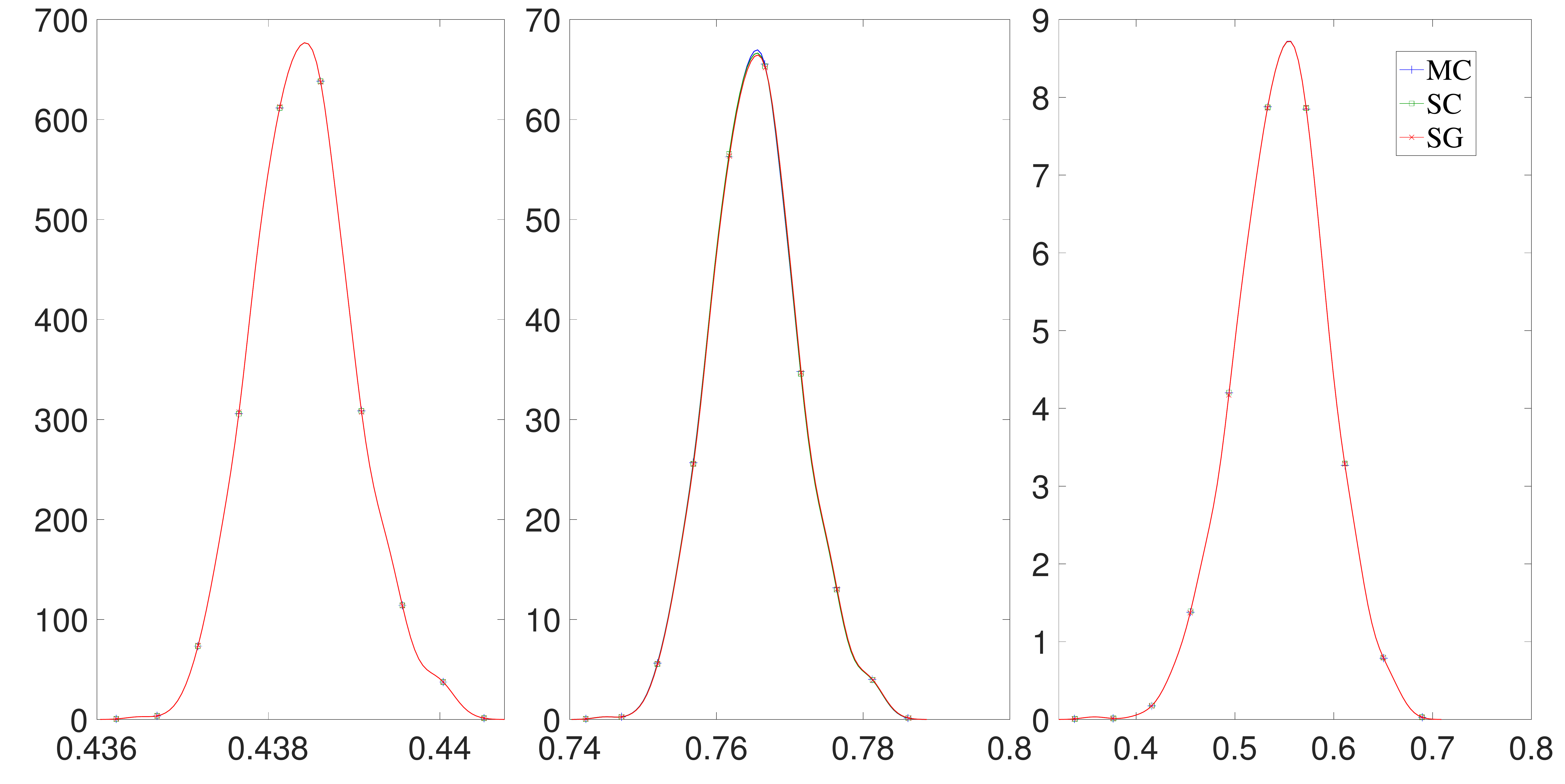}\\
\includegraphics[width=6.5cm]{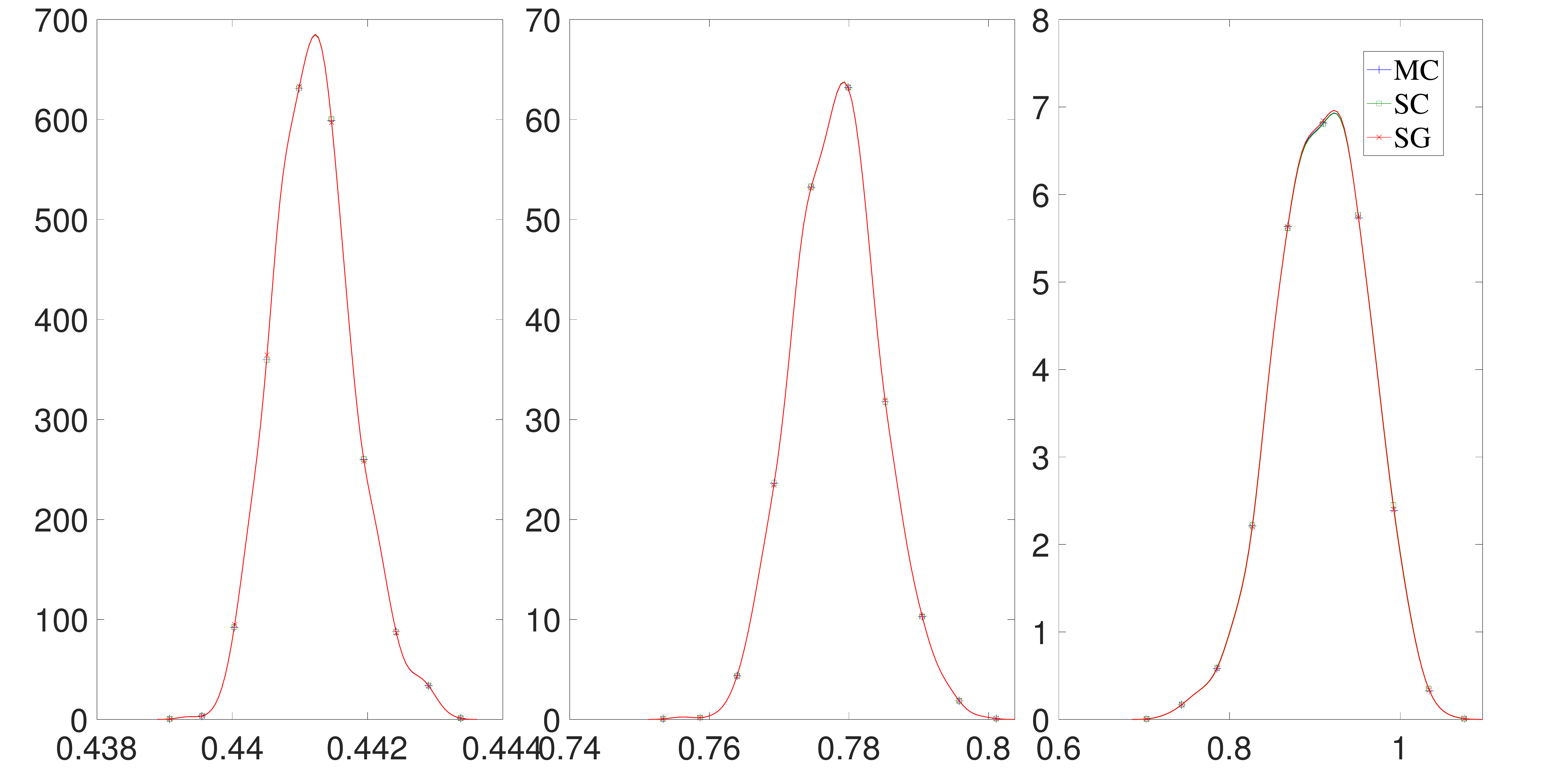} &
\includegraphics[width=6.5cm]{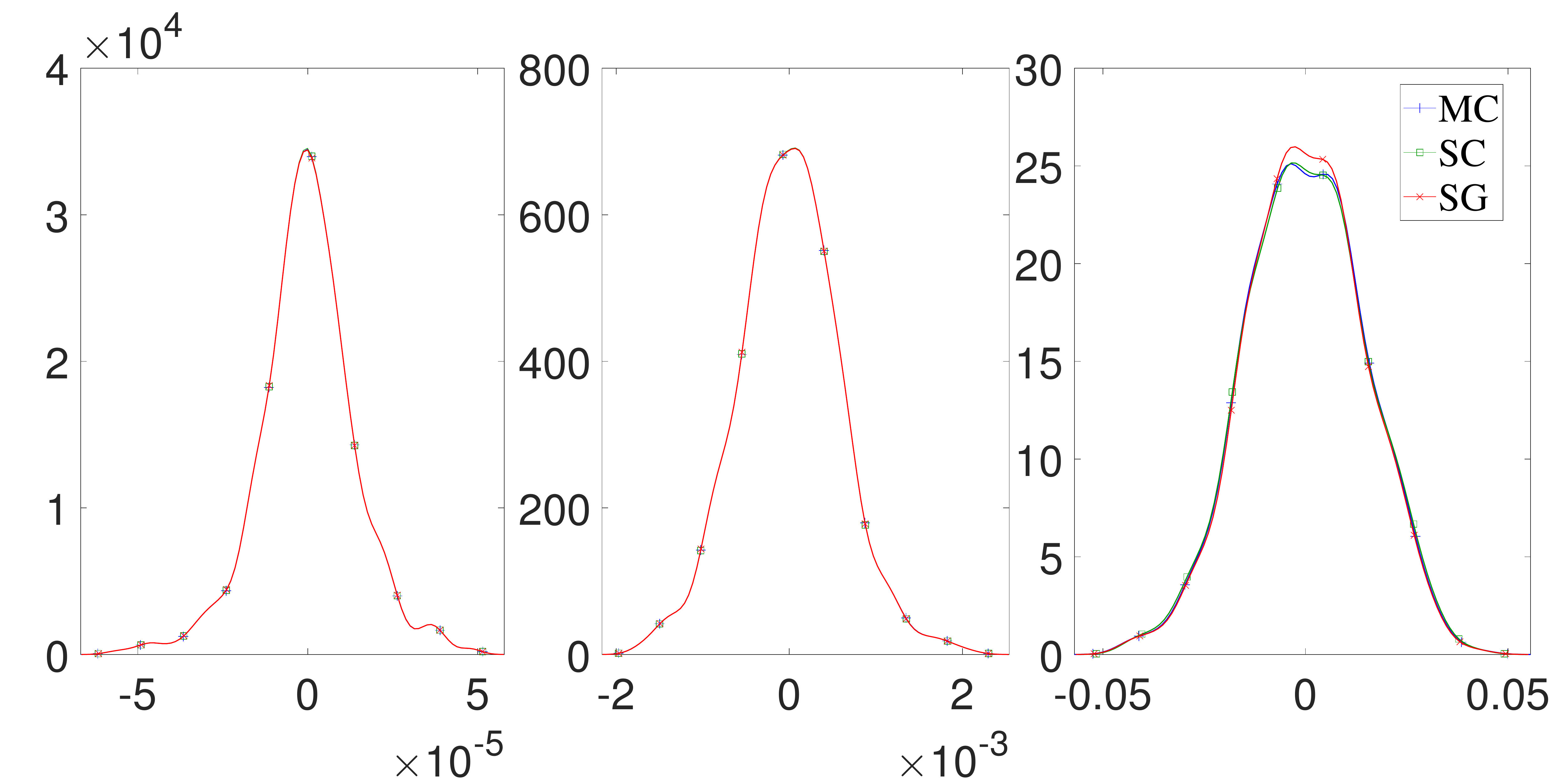}
\end{tabular}
\end{center}
\caption{Estimated probability density functions at times $0.1$s, $1$s and
$10$s (left to right, all panels) of the horizontal velocity at points with
coordinates $(4.0100,-0.4339)$ (top left), $(4.0100,0.4339)$ (bottom left),
and of the horizontal (top right) and vertical (bottom right) velocities at
the point $(3.6436,0)$ for mean Reynolds number $\operatorname*{Re}_{1}=100$
and $CoV=10\%$.}%
\label{fig:QoI_100_10}%
\end{figure}

\begin{figure}[ptbh]
\begin{center}`
\includegraphics[width=12cm]{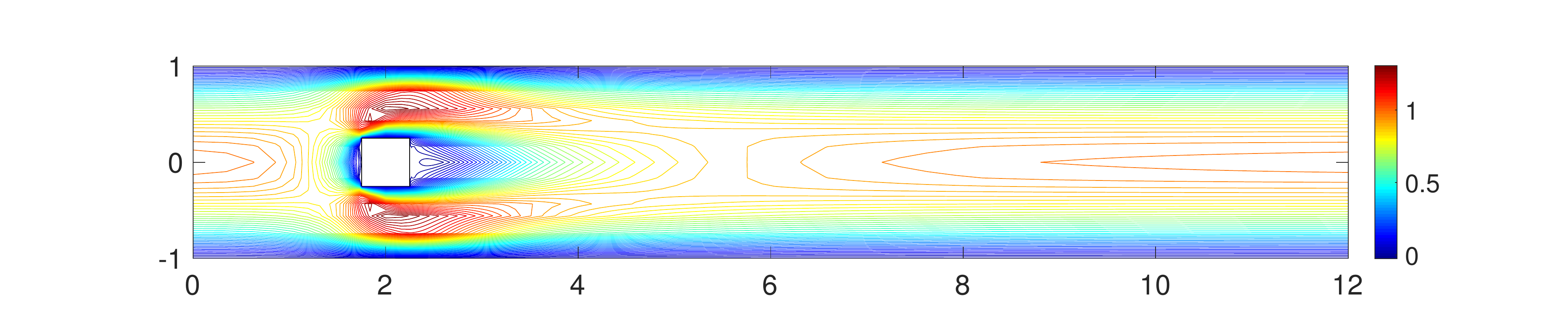}
\includegraphics[width=12cm]{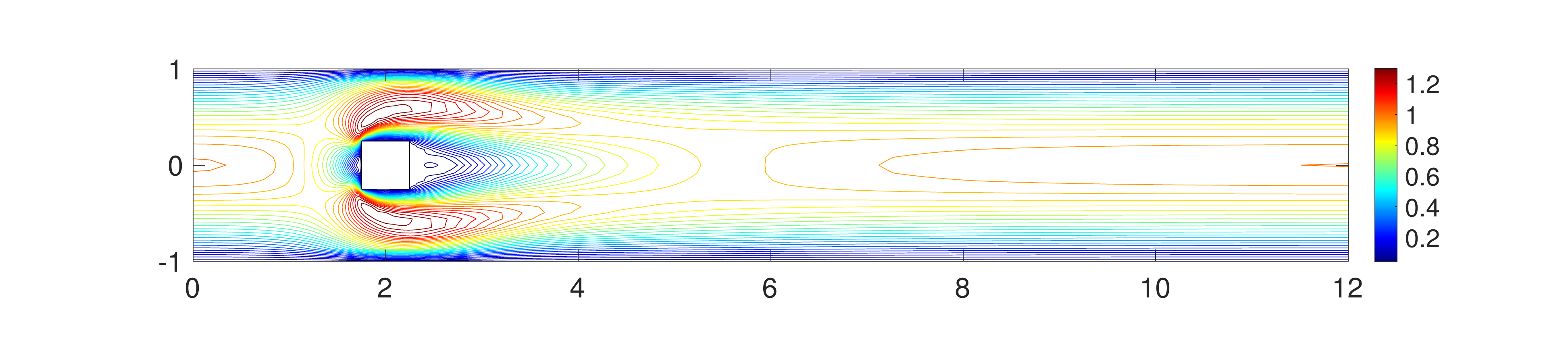}
\end{center}
\caption{Mean horizontal velocity obtained using the stochastic Galerkin methods for the steady-state problem (top),
and at time $100$s (bottom) with mean Reynolds number $\operatorname*{Re}_{1}=100$ and $CoV=10\%$.}%
\label{fig:mean_ux_100_10_steady}%
\end{figure}

\begin{figure}[ptbh]
\begin{center}
\includegraphics[width=12cm]{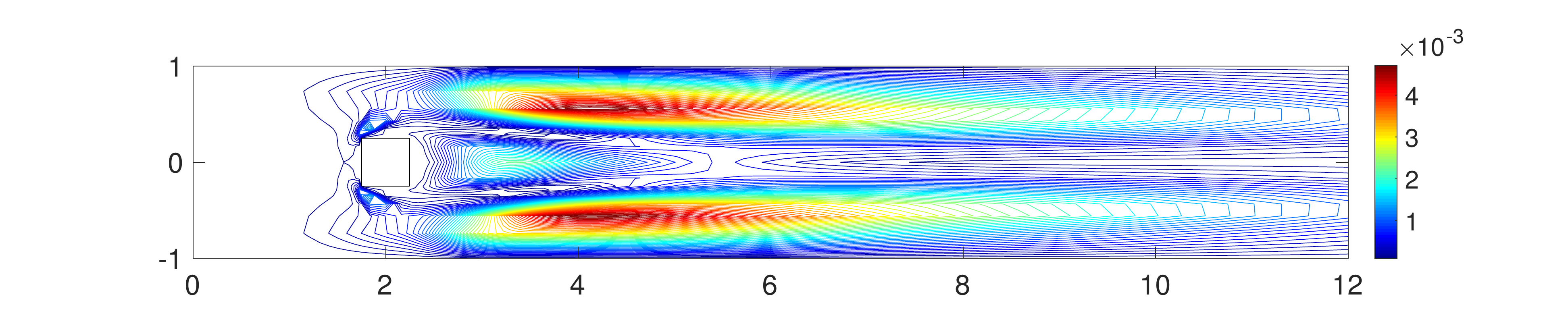}
\includegraphics[width=12cm]{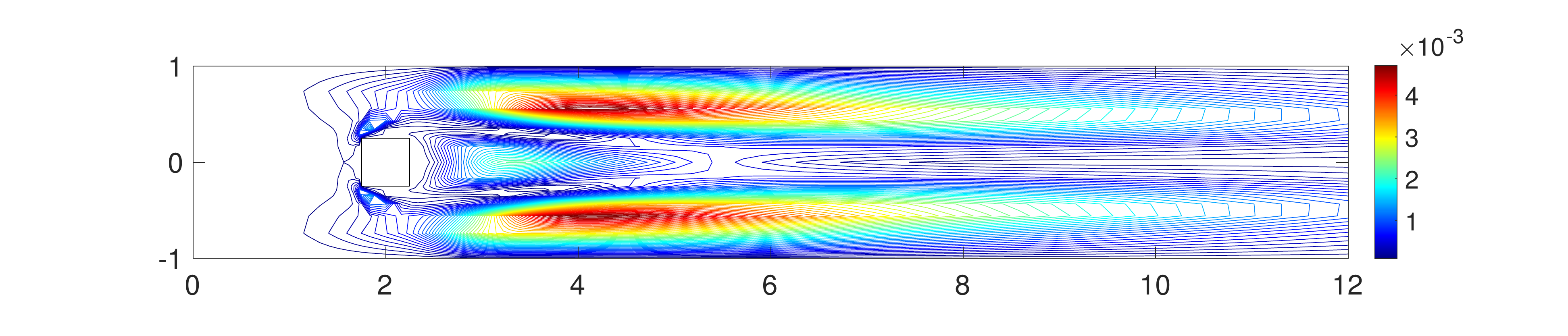}
\end{center}
\caption{Variance of the horizontal velocity obtained using the stochastic Galerkin methods for the steady-state problem (top),
and at time $100$s (bottom) with mean Reynolds number $\operatorname*{Re}_{1}=100$ and $CoV=10\%$.}%
\label{fig:variance_ux_100_10_steady}%
\end{figure}

Now, let us consider first the case of $\operatorname*{Re}_{1}=100$ and
$CoV=10\%$.
Figure~\ref{fig:gPC_100_10} shows the evolution of the gPC\ coefficients of
the horizontal velocity, and the symbols $\square$ and $\times$ represent the
results of Monte Carlo and stochastic collocation at some of the time
barriers. It can be seen that all methods are in agreement.
Figure~\ref{fig:mean_velocity_100_10} shows the mean horizontal velocity,
Figure~\ref{fig:variance_ux_100_10} the variance of the horizontal velocity,
and Figure~\ref{fig:variance_uy_100_10} the variance of the vertical velocity
at times $0.1$s, $1$s and $10$s. From these figures it can be seen that the
flow quickly evolves during the first second, and the later changes are
relatively less dramatic. It can be seen that there is symmetry in all the
quantities, the mean values are essentially the same as we would expect in the
deterministic case~\cite{Elman-2014-FEF}, and the variance of the horizontal
velocity component is evolving to be concentrated in two \textquotedblleft
eddies" and it is larger than the variance of the vertical velocity component.
\textcolor{black}{In fact, it appears that all quantities are already at time $10$s close to the steady state, 
see also Figures~\ref{fig:mean_ux_100_10_steady} and~\ref{fig:variance_ux_100_10_steady}.}
A different perspective on the
solution is given by Figure~\ref{fig:QoI_100_10}, which displays evolution of
the probability density function (pdf) estimates in several points of the
domain at times $0.1$s, $1$s and $10$s. The left panels show the pdf estimates
of the velocity in $x$ direction at points with coordinates $(4.0100,-0.4339)$
(top), $(4.0100,0.4339)$ (bottom), where the variance of the velocity is
relatively large cf. Figure~\ref{fig:variance_ux_100_10}. The right panels
show the estimates at point $(3.6436,0)$ which is slightly downstream from the
obstacle: the estimate in the $x$ direction in the top panel and the estimate
in the $y$ direction in the bottom panel. The results were obtained using
\textsc{Matlab}'s \texttt{ksdensity} function. It can be seen that the changes
of the mean values of the pdf estimates are relatively large during the first
second, and then the uncertainty gradually increases and the supports of the
pdf estimates grow as the solution evolves away from the deterministic initial
condition and the effect of the stochastic viscosity becomes evident.


\begin{figure}[ptbh]
\begin{center}
\includegraphics[width=10cm]{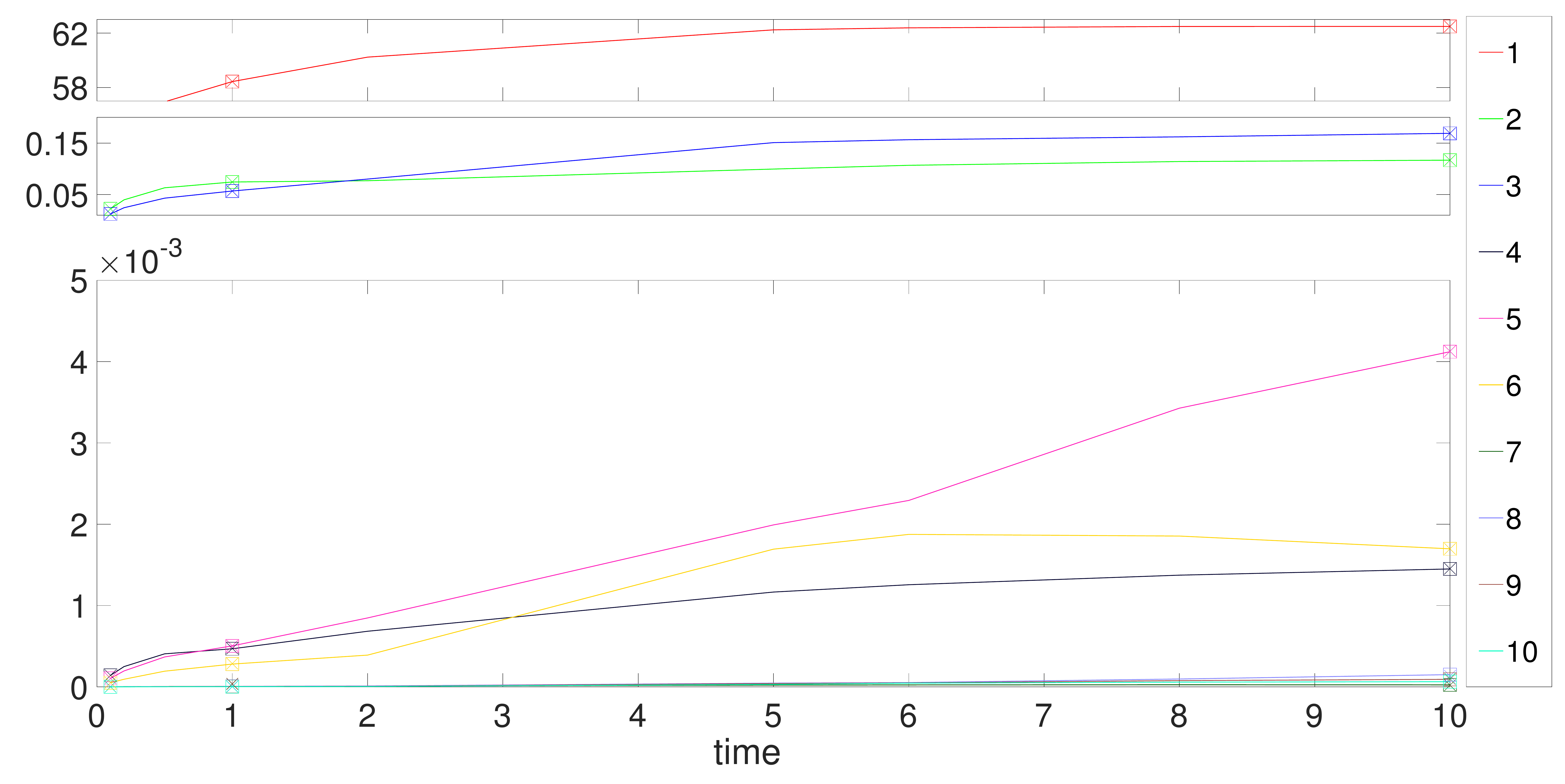}
\end{center}
\caption{Evolution of the gPC coefficients corresponding to the horizontal
velocity in terms of $\ell_{2}$-norm for mean Reynolds number
$\operatorname*{Re}_{1}=300$ and $CoV=1\%$. The symbols $\square$ and $\times$
represent the results of the Monte Carlo and stochastic collocation,
respectively, at times $0.1$s, $1$s and $10$s.}%
\label{fig:gPC_300_1}%
\end{figure}

\begin{figure}[ptbh]
\begin{center}
\includegraphics[width=12cm]{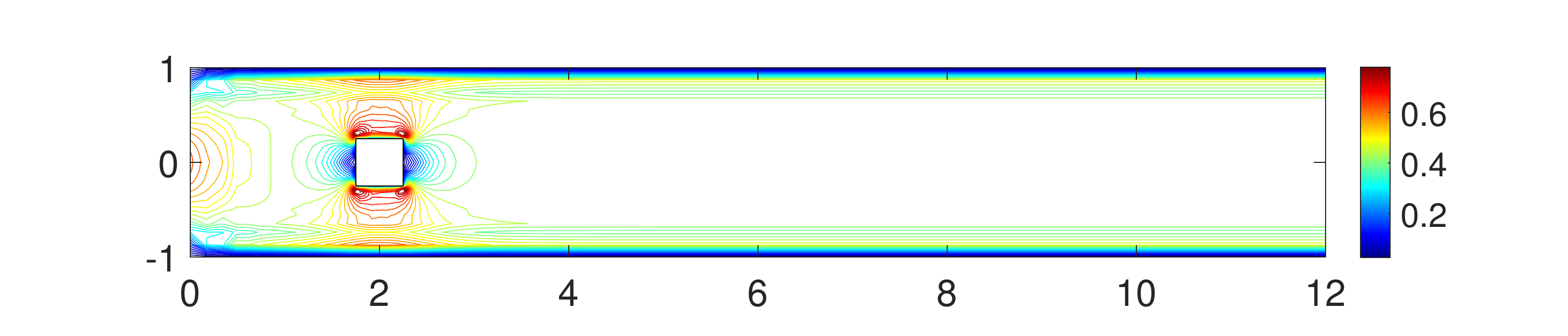}
\includegraphics[width=12cm]{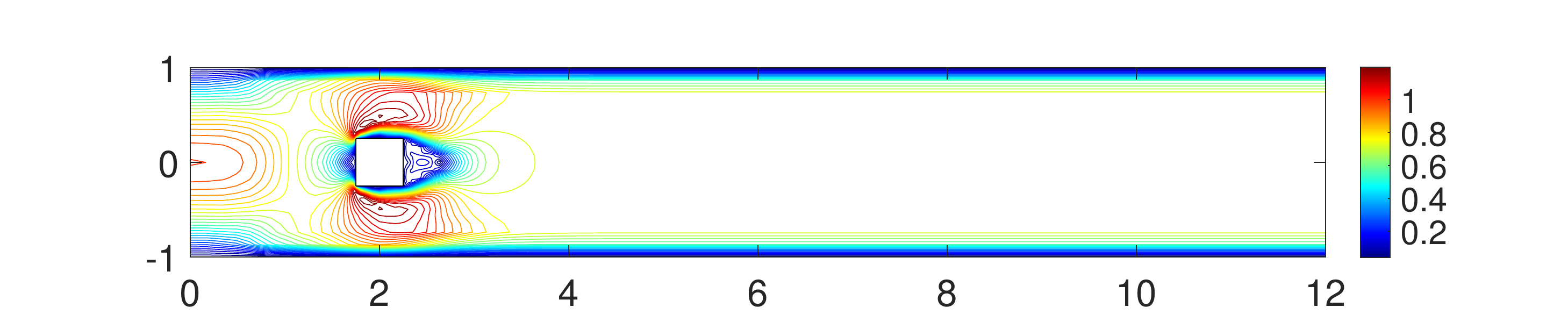}
\includegraphics[width=12cm]{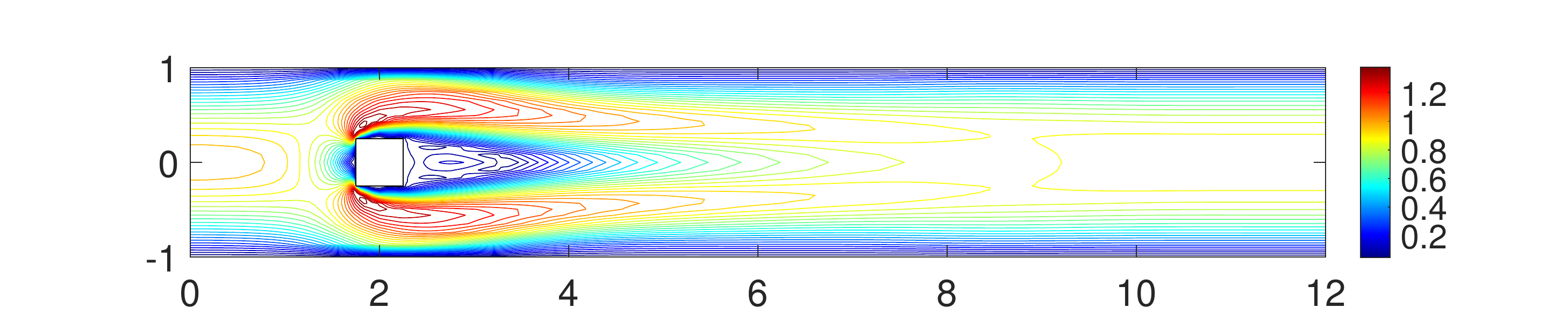}
\end{center}
\caption{Mean horizontal velocity at times $0.1$s (top), $1$s (center) and
$10$s (bottom) for mean Reynolds number $\operatorname*{Re}_{1}=300$ and
$CoV=1\%$.}%
\label{fig:mean_velocity_300_1}%
\end{figure}

\begin{figure}[ptbh]
\begin{center}
\includegraphics[width=12cm]{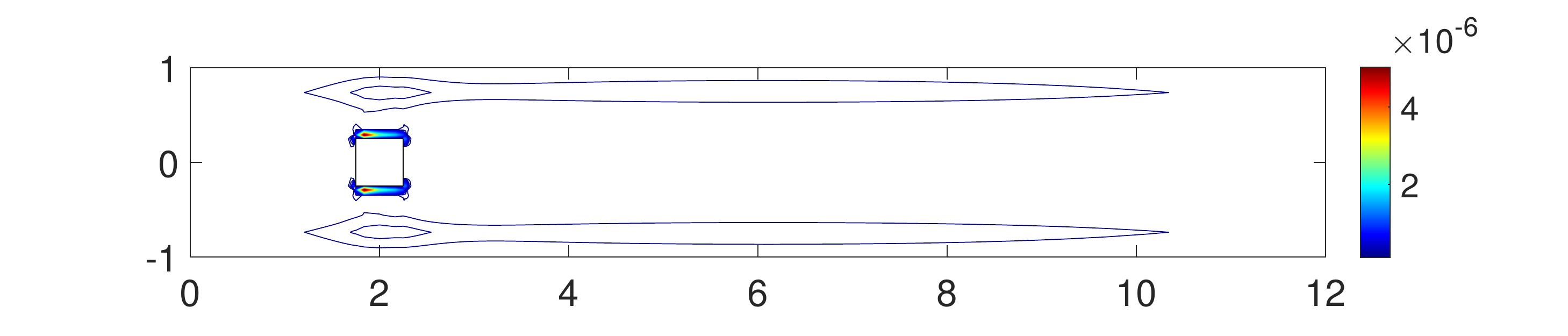}
\includegraphics[width=12cm]{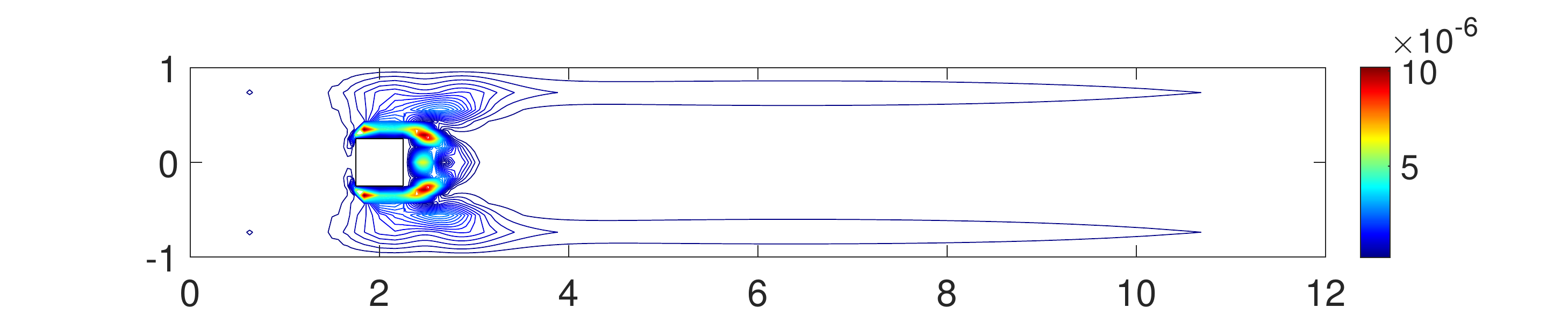}
\includegraphics[width=12cm]{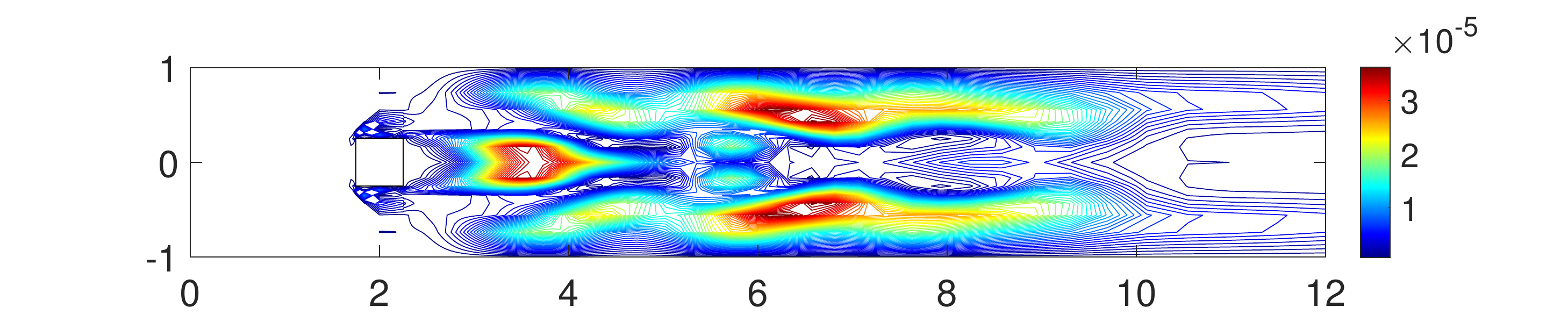}
\end{center}
\caption{Variance of the horizontal velocity at times $0.1$s (top), $1$s
(center) and $10$s (bottom) for mean Reynolds number $\operatorname*{Re}%
_{1}=300$ and $CoV=1\%$.}%
\label{fig:variance_ux_300_1}%
\end{figure}

\begin{figure}[ptbh]
\begin{center}
\includegraphics[width=12cm]{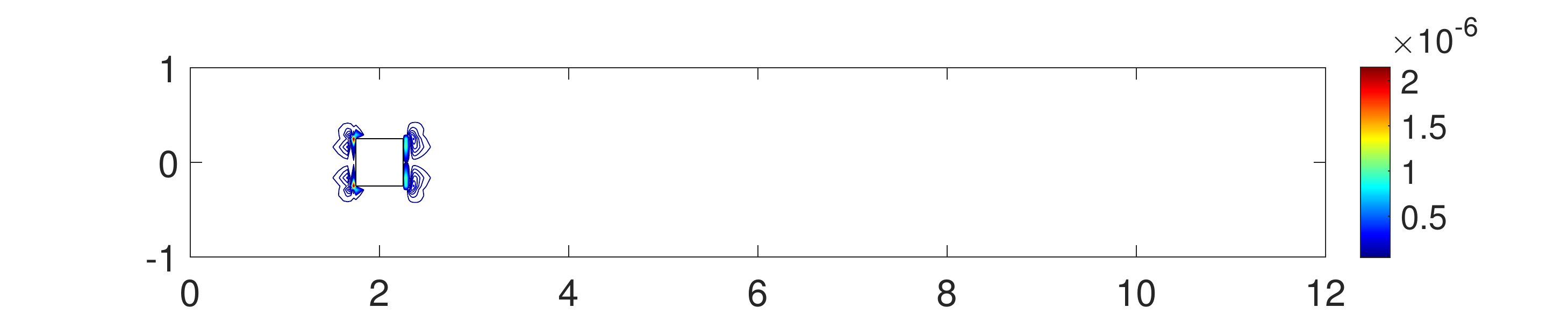}
\includegraphics[width=12cm]{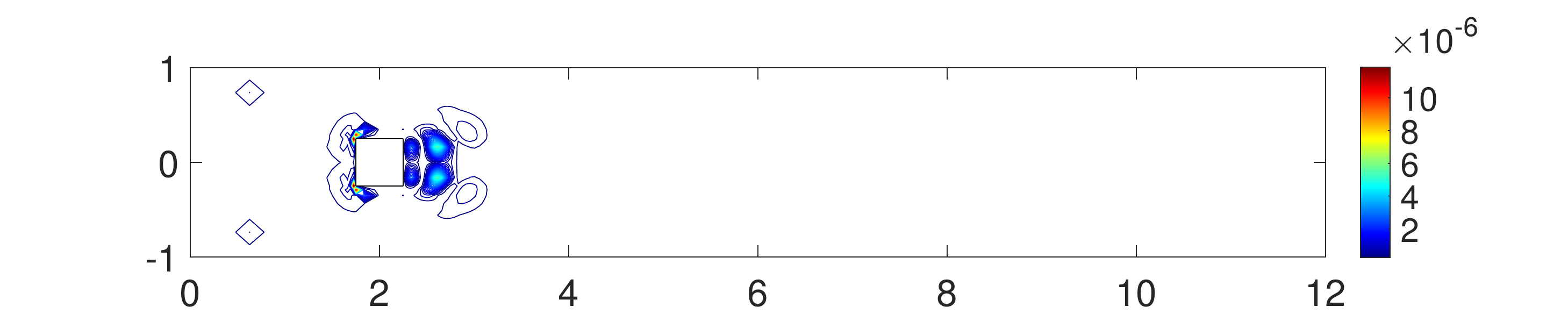}
\includegraphics[width=12cm]{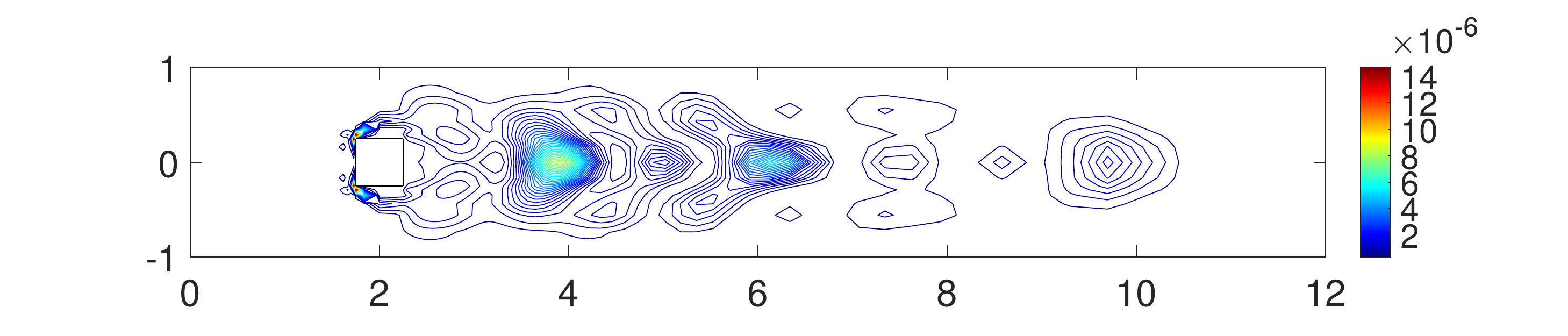}
\end{center}
\caption{Variance of the vertical velocity at times $0.1$s (top), $1$s
(center) and $10$s (bottom) for mean Reynolds number $\operatorname*{Re}%
_{1}=300$ and $CoV=1\%$.}%
\label{fig:variance_uy_300_1}%
\end{figure}

\begin{figure}[ptbh]
\begin{center}%
\begin{tabular}
[c]{c@{\hskip-0.2cm}c}%
\includegraphics[width=6.5cm]{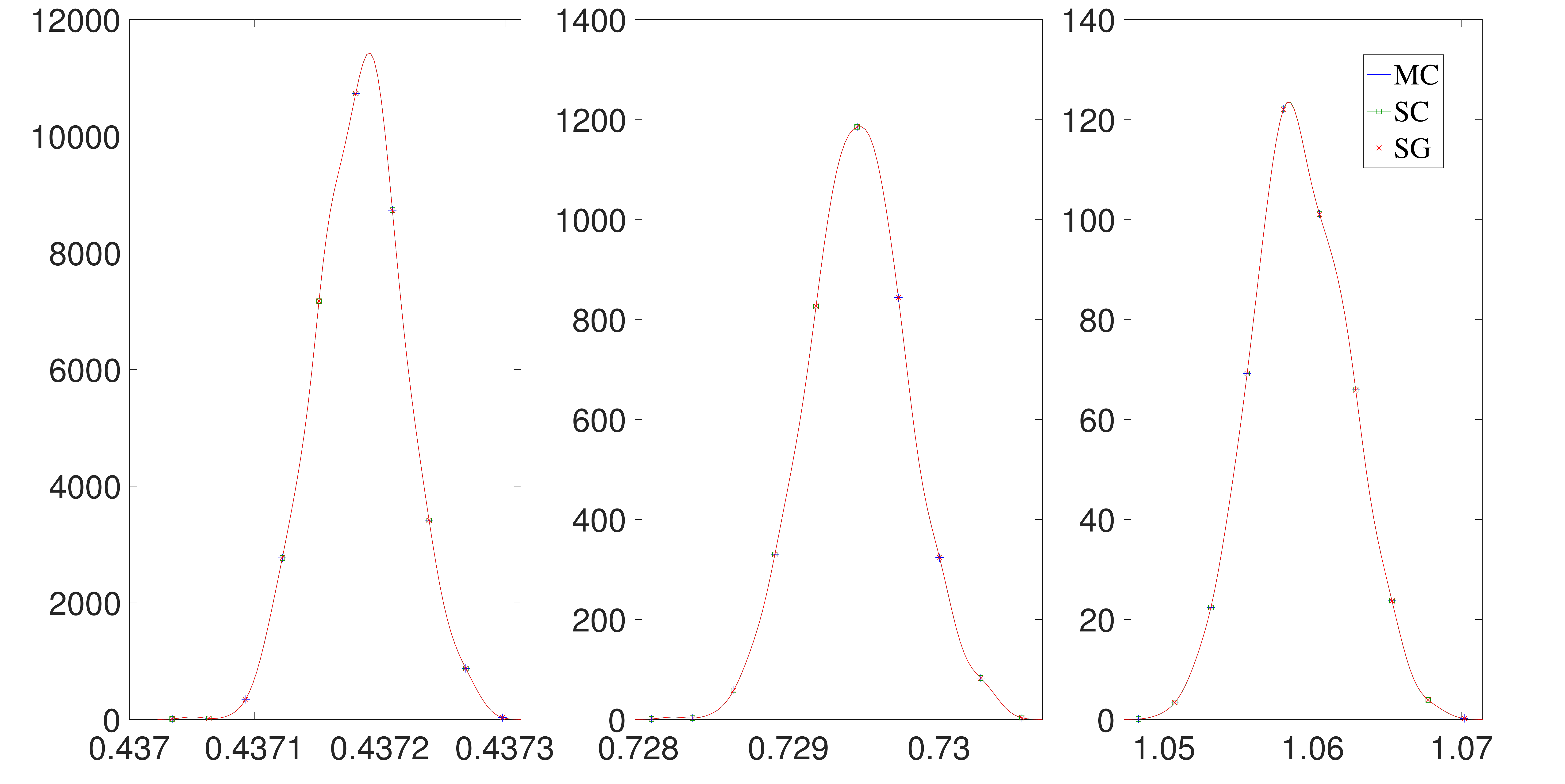} &
\includegraphics[width=6.5cm]{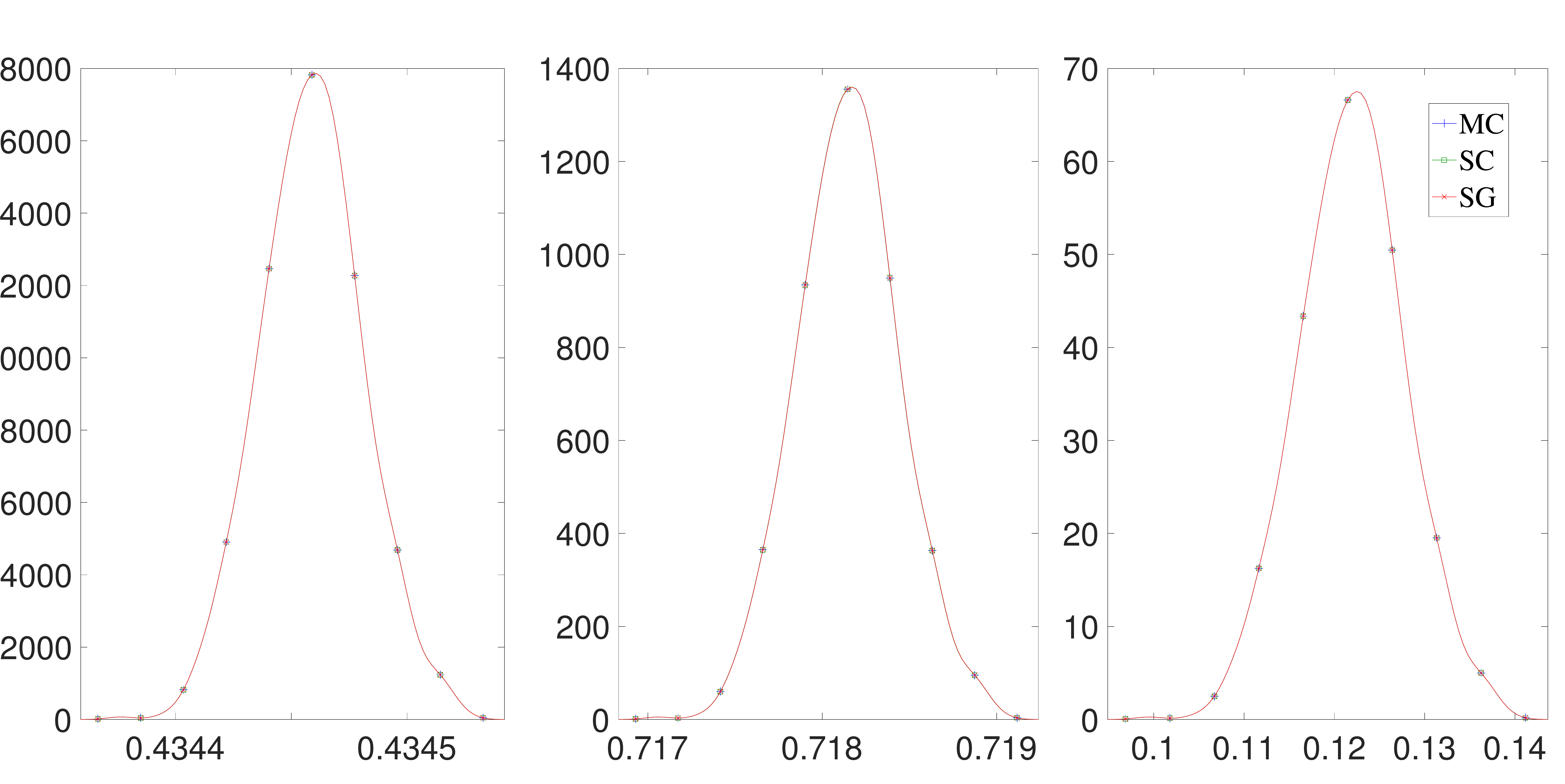}\\
\includegraphics[width=6.5cm]{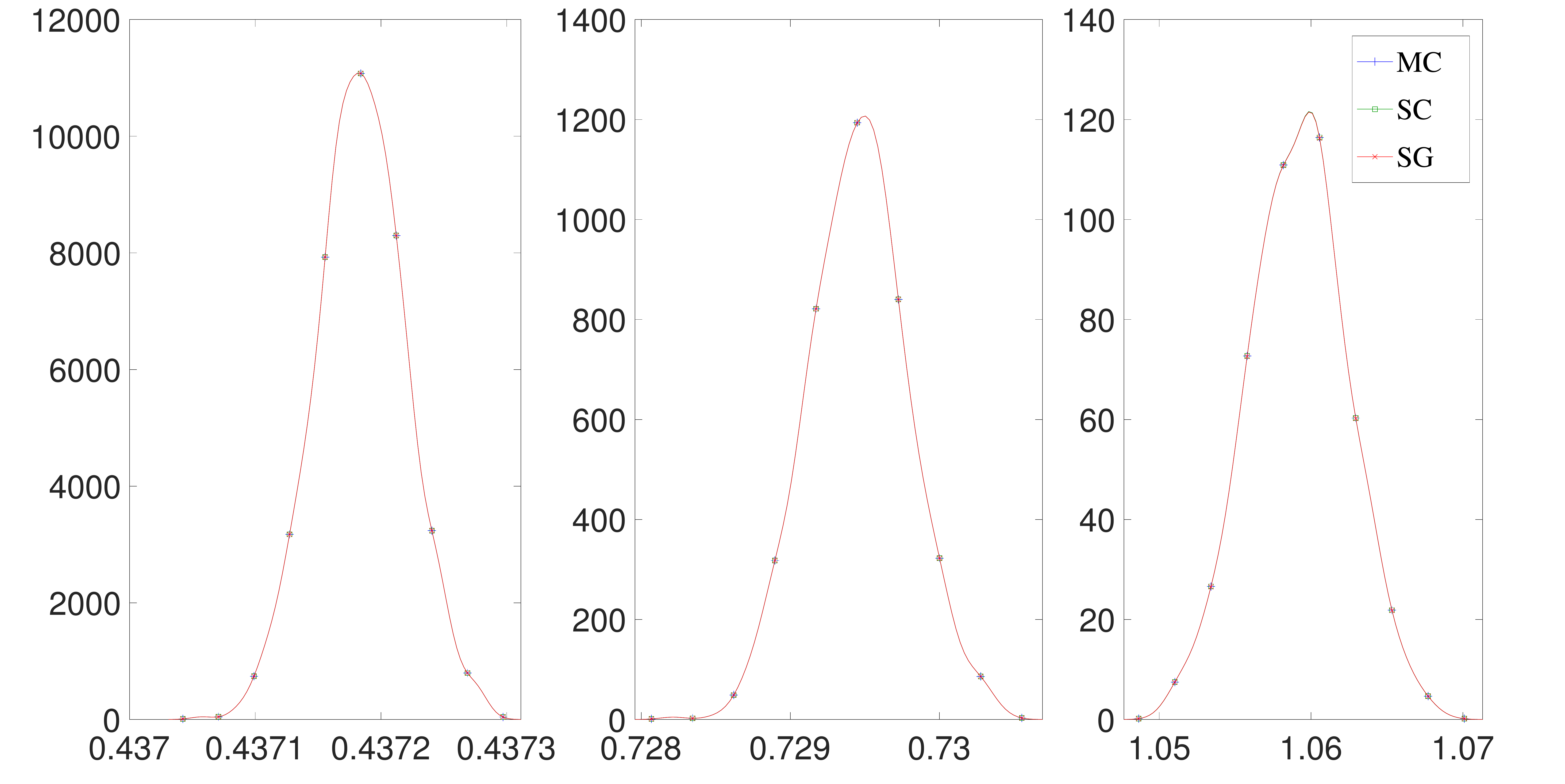} &
\includegraphics[width=6.5cm]{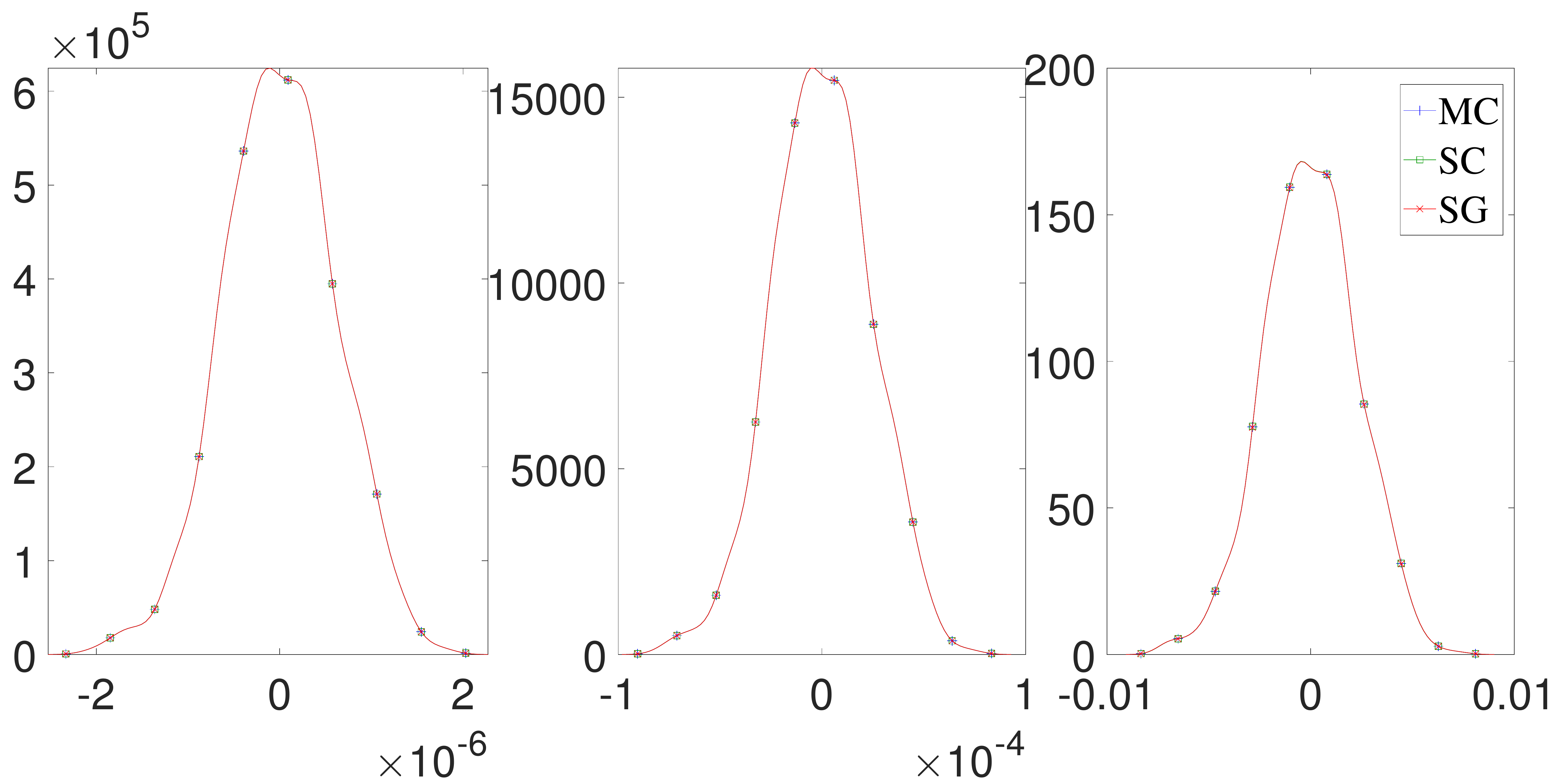}
\end{tabular}
\end{center}
\caption{Estimated probability density functions at times $0.1$s, $1$s and
$10$s (left to right, all panels) of the horizontal velocity at points with
coordinates $(4.0100,-0.4339)$ (top left), $(4.0100,0.4339)$ (bottom left),
and of the horizontal (top right) and vertical (bottom right) velocities at
the point $(3.6436,0)$ for mean Reynolds number $\operatorname*{Re}_{1}=300$
and $CoV=1\%$.}%
\label{fig:QoI_300_1}%
\end{figure}

Next, let us consider the case of $\operatorname*{Re}_{1}=300$ and $CoV=1\%$.
Figure~\ref{fig:gPC_300_1} shows the evolution of the gPC coefficients of the
horizontal velocity. It can be seen that with increased $\operatorname*{Re}%
_{1}$ it takes more time for the flow to develop, including the stochastic
components of the solution despite lower $CoV$ than in the previous problem.
Again, all methods are in agreement. Figure~\ref{fig:mean_velocity_300_1} then
shows the mean horizontal velocity, Figure~\ref{fig:variance_ux_300_1} the
variance of the horizontal velocity, and Figure~\ref{fig:variance_uy_300_1}
the variance of the vertical velocity, all at times $0.1$s, $1$s and $10$s.
The mean quantities are quite similar to what would be expected in the
deterministic case, and the variances reflect on more complex behavior of the
fluid at the higher value of $\operatorname*{Re}_{1}$. Finally,
Figure~\ref{fig:QoI_300_1} displays evolution of the probability density
function (pdf) estimates at the same set of points of the domain and times as
Figure~\ref{fig:QoI_100_10}, and all three methods are again in agreement.

\textcolor{black}{Finally, we compare the results of the stochastic Galerkin method applied to the steady-state problem 
with mean Reynolds number $\operatorname*{Re}_{1}=100$ and $CoV=10\%$, which was 
studied by Soused\'{i}k and Elman in~\cite{Sousedik-2016-SGM}, 
and the results of the long-term integration at time~$100$s. 
Specifically, a comparison of the mean horizontal velocity is shown in Figure~\ref{fig:mean_ux_100_10_steady}, and 
Figure~\ref{fig:variance_ux_100_10_steady} displays the variance of the horizontal velocity.
By comparing the two figures, it can be seen that the results are virtually identical.}

\subsection{Preconditioning of the Oseen problem}

\label{sec:preconditioner}
\begin{figure}[ptbh]
\begin{center}
\includegraphics[width=12cm]{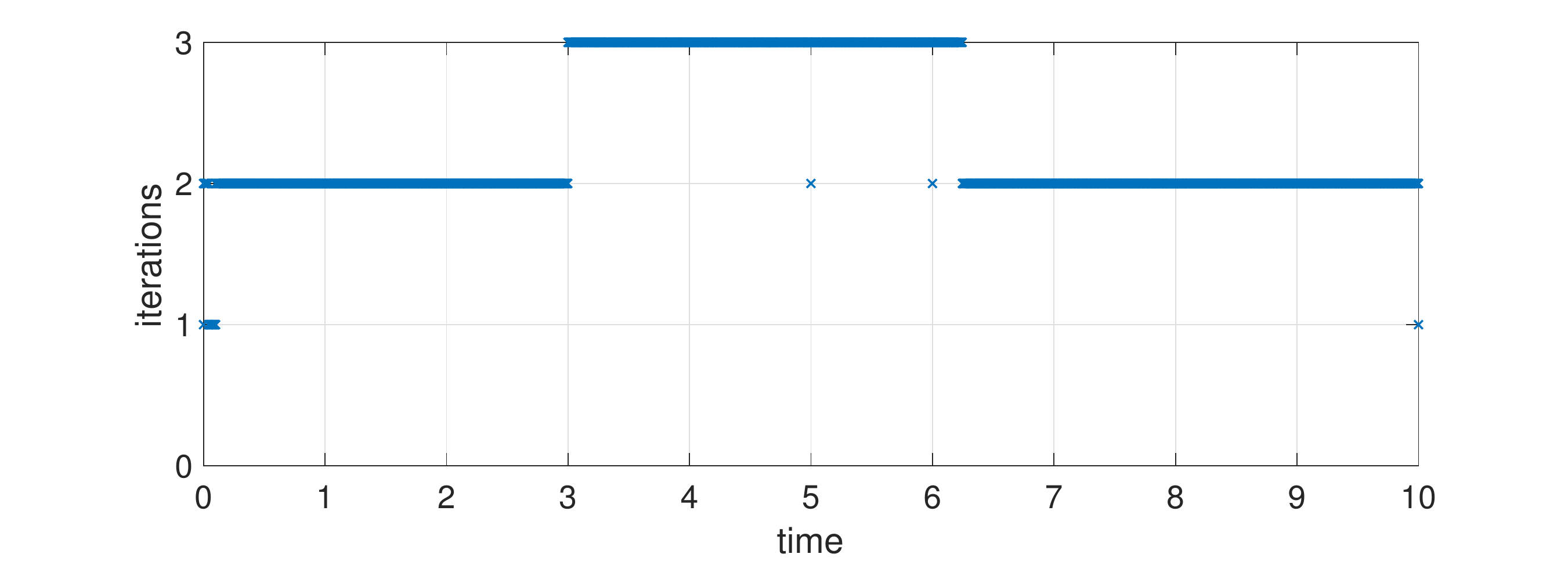} \\
\includegraphics[width=12cm]{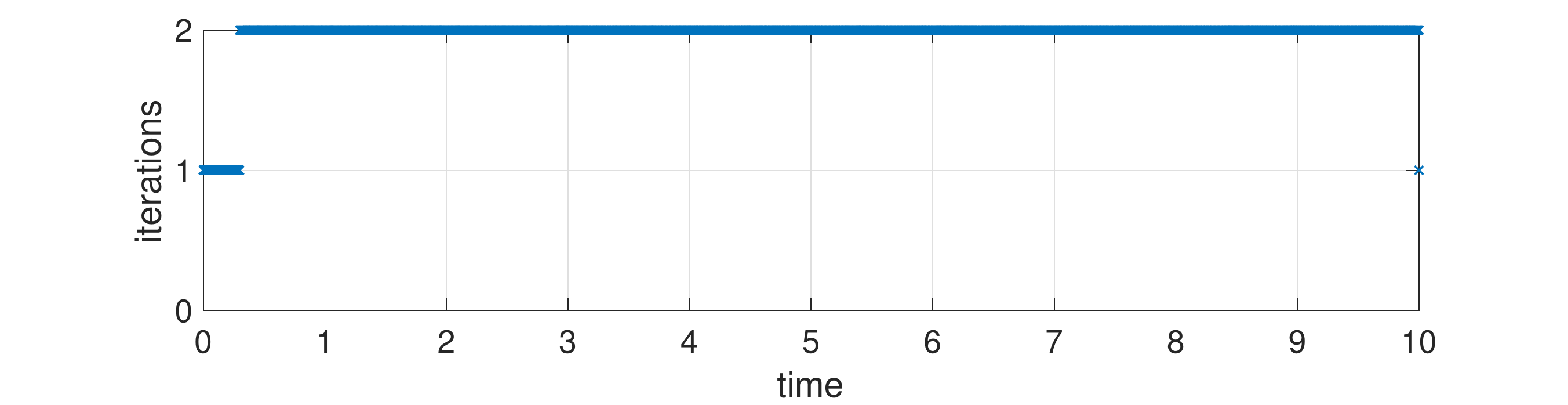}
\end{center}
\caption{Numbers of fGMRES iterations with the exact mean-based preconditioner for
mean Reynolds number $\operatorname*{Re}_{1}=100$, $CoV=10\%$ (top),
and $\operatorname*{Re}_{1}=300$, $CoV=1\%$ (bottom).}%
\label{fig:gmres_it_LU}%
\end{figure}

\begin{figure}[ptbh]
\begin{center}
\includegraphics[width=12cm]{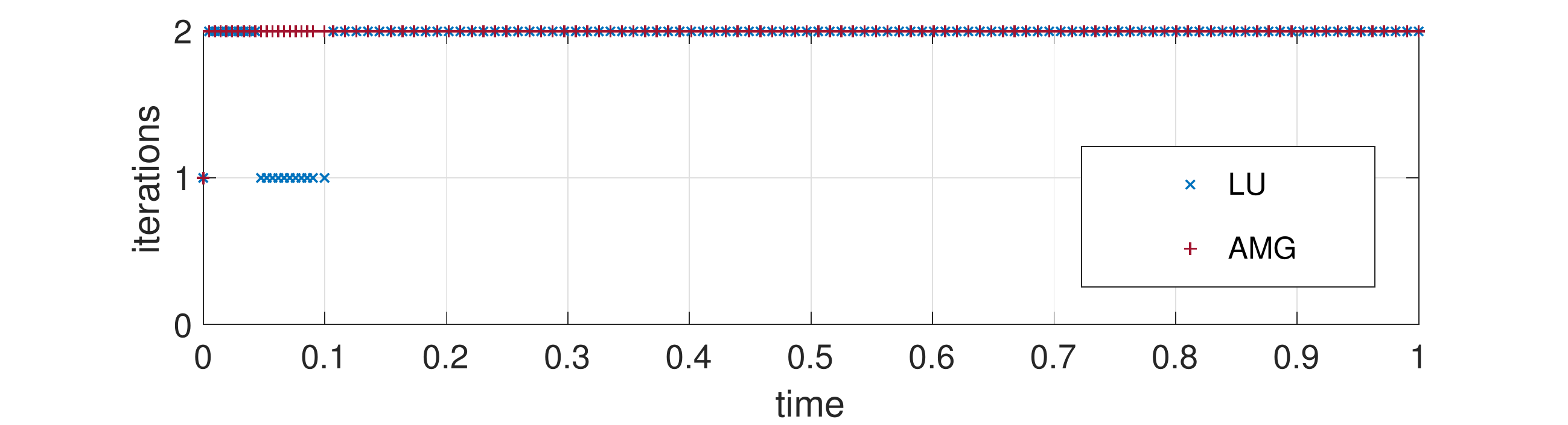} \\
\includegraphics[width=12cm]{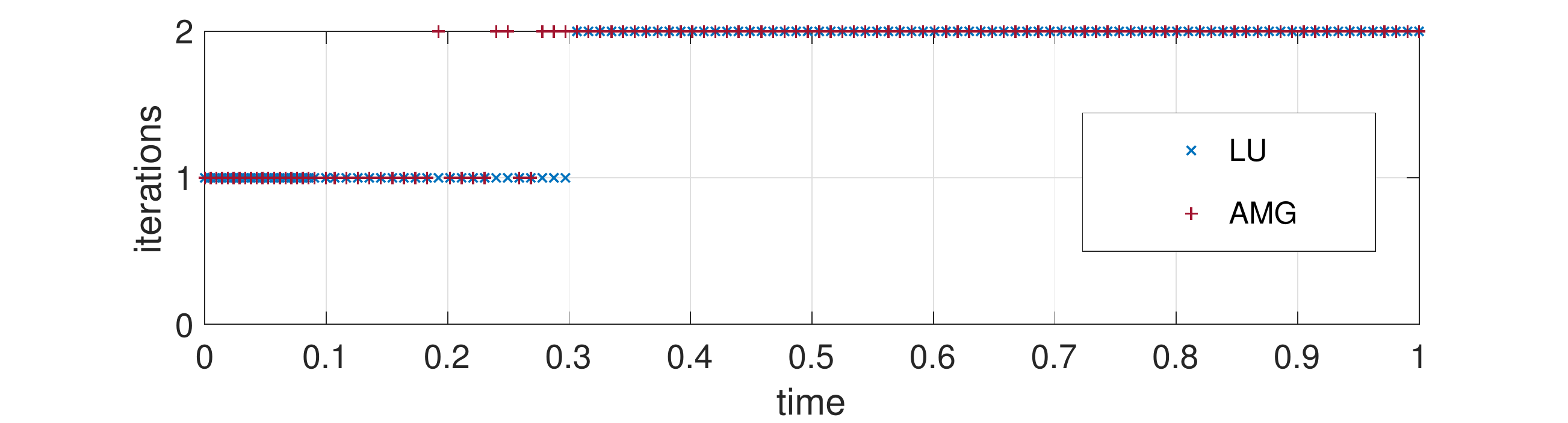} 
\end{center}
\caption{A comparison of the exact mean-based preconditioner (LU) and its PCD iterated variant (AMG)
in terms of the numbers of fGMRES iterations for 
mean Reynolds number $\operatorname*{Re}_{1}=100$, $CoV=10\%$ (top),
and $\operatorname*{Re}_{1}=300$, $CoV=1\%$ (bottom).}%
\label{fig:gmres_it_LU-AMG}%
\end{figure}

The solution of the Oseen problem~(\ref{eq:Oseen-s})
in each time step of the stochastic Galerkin method is a computationally
expensive task. Therefore, use of a preconditioned Krylov subspace method may
be preferred over a direct solver. To this end, we used the right-preconditioned
flexible GMRES (fGMRES) method~\cite{Saad-1993-FIP} with the so-called
mean-based preconditioner $\mathcal{M}_{1}^{-1}:\mathbf{R}\longmapsto
\mathbf{V}$, which entails solving 
a linear system
\begin{equation}
\mathcal{M}_{1}\mathbf{V}=\mathbf{R,}\label{eq:alg-MB1}%
\end{equation}
where $\mathbf{R}$ and $\mathbf{V}$ are the matricized coefficients of the
gPC\ expansions, cf.~(\ref{eq:U}). Specifically, $\mathcal{M}_{1}^{-1}$
denotes an action of the pressure convection-diffusion (PCD)\ preconditioner,
see~\cite[Section~$3$]{Kay-2010-ATS} and~\cite[Section 9.2.2]{Elman-2014-FEF},
which is motivated by the block inverse of the matrix$~\mathcal{F}_{1}^{n+1}$
in~(\ref{eq:stochastic-saddle-point}). It can be specifically written as
\begin{equation}
\mathcal{M}_{1}^{-1}=\left[
\begin{array}
[c]{cc}%
\left(  \mathbf{F}_{1}^{n+1}\right)  ^{-1} & \left(  \mathbf{F}_{1}%
^{n+1}\right)  ^{-1}\mathbf{B}^{T}\mathbf{X}^{-1}\\
\mathbf{0} & \mathbf{-X}^{-1}%
\end{array}
\right]  ,\label{eq:M_inv}%
\end{equation}
where $\mathbf{F}_{1}^{n+1}$ is the matrix from~(\ref{eq:F_1-s}), and
$\mathbf{X}^{-1}$ is the pressure convection-diffusion~term
\[
\mathbf{X}^{-1}\mathbf{=A}_{p}^{-1}\mathbf{F}_{p}^{n+1}\mathbf{M}_{p}^{-1}.
\]
First, we used LU\ factorizations of the matrices from~(\ref{eq:M_inv}), which
are updated in each time step. Since the solves with the (mean)
matrix~$\mathcal{M}_{1}$ are thus exact, this illustrates the approximation
properties of the mean-based preconditioner. Then, we also used the
IFISS\ implementation of the PCD\ iterated preconditioner, in which the solves
involving both $\mathbf{F}_{1}^{n+1}$ and $\mathbf{A}_{p}=\mathbf{BT}%
^{-1}\mathbf{B}^{T}$, where $\mathbf{T}$ is the diagonal of the velocity mass
matrix, are replaced by a single V-cycle of AMG\ using the IFISS\ default
parameters, and the solve with the pressure matrix$~\mathbf{M}_{p}$ is
effected by five Chebyshev iterations, see~\cite[Section~10.3]{Elman-2014-FEF}%
. The construction of the matrix $\mathbf{F}_{p}^{n+1}$ is described
in~\cite[Chapter~9]{Elman-2014-FEF}. We note that the AMG\ implementation is
based on \texttt{HSL\_MI20}~\cite{Boyle-2010-HSL}. All tests started with a
zero initial iterate and stopped when the relative residual was reduced to
$10^{-8}$ in the Euclidean norm. 
The numbers of fGMRES iterations for solves
in the time interval $\left[  0,10s\right] $ with the exact mean-based
preconditioner~(\ref{eq:alg-MB1}) are shown in Figure~\ref{fig:gmres_it_LU}.
It can be seen that at most three iterations were needed in all steps.  
A comparison of the exact mean-based preconditioner (LU) and its PCD iterated variant (AMG) is illustrated by 
Figure~\ref{fig:gmres_it_LU-AMG}. 
It can be seen that the numbers of iterations are the same in most cases, or 
it takes at most one extra step for the PCD iterated variant to converge. 
Thus both exact and iterated
versions of the mean-based preconditioner are suitable for the
problems studied in our numerical experiments.

\section{Conclusion}

\label{sec:conclusion} 
\textcolor{black}{We studied the time-dependent Navier\textendash Stokes
equations with stochastic viscosity, which was given in terms of a polynomial chaos expansion. 
For this problem, we developed a stochastic Galerkin method with
adaptive, mean-informed time stepping. 
We applied the method to a popular benchmark problem given by a flow around an obstacle,
and we compared the solution of the time-dependent problem after the transient to that of the corresponding steady-state problem.
Next, since the time-stepping scheme is fully implicit, a linear solve with the stochastic Galerkin matrix is required in each time step. 
Use of direct solvers may be prohibitive due to the large size of the systems, 
and in fact it is even not desirable to form the matrices explicitly. 
Therefore, we also formulated a preconditioner, which is used by the right-preconditioned flexible GMRES method,
and allows to solve the stochastic Galerkin systems efficiently. 
We studied two variants of the preconditioner. 
The first variant is based on exact factorization of the matrix corresponding to the underlying mean problem,
and the second one was an iterated variant by means of an algebraic multigrid solver. 
In the numerical experiments we observed that the performance of the exact and iterated variants of the preconditioner was virtually identical, 
and only a couple of GMRES iterations were needed for convergence in all time steps. 
Therefore, the proposed stochastic Galerkin method is designed as a wrapper around an existing code for the corresponding deterministic problem, 
and in fact an efficient solver for the deterministic problem is the essential component also for the method presented in this study. 
Finally, we also compared the stochastic Galerkin solution with the
stochastic collocation and Monte Carlo solutions, and 
we observed an excellent agreement for all problems studied in our numerical experiments. 
}


\bibliographystyle{plain}
\bibliography{gPC_time.bib}

\end{document}